\documentclass[12pt]{amsart}
\usepackage{amssymb}
\usepackage{mathrsfs}
\usepackage[english]{babel}
\input xy
\xyoption{all}

\begin{document}

\newtheoremstyle{mytheorem}
  {}
  {}
  {\slshape}
  {}
  {\scshape}
  {.}
  { }
  {}

\newtheoremstyle{mytheoremintro}
  {}
  {}
  {\slshape}
  {}
  {\scshape}
  {.}
  { }
  {}

\newtheoremstyle{mydefinition}
  {}
  {}
  {\upshape}
  {}
  {\scshape}
  {.}
  { }
  {}

\theoremstyle{mytheoremintro}
\newtheorem{thmintro}{Theorem}
\newtheorem{corintro}[thmintro]{Corollary}
\theoremstyle{mytheorem}
\newtheorem{lemma}{Lemma}[section]
\newtheorem{prop}[lemma]{Proposition}
\newtheorem{cor}[lemma]{Corollary}
\newtheorem{thm}[lemma]{Theorem}
\theoremstyle{mydefinition}
\newtheorem{rem}[lemma]{Remark}
\newtheorem*{rem*}{Remark}
\newtheorem{conj}[lemma]{Conjecture}
\newtheorem{rems}[lemma]{Remarks}
\newtheorem{defi}[lemma]{Definition}
\newtheorem*{defi*}{Definition}
\newtheorem{defis}[lemma]{Definitions}
\newtheorem{exo}[lemma]{Example}
\newtheorem{exos}[lemma]{Examples}
\renewcommand{\thethmintro}{\Alph{thmintro}}
\numberwithin{equation}{section}

\newcommand{\bibURL}[1]{{\unskip\nobreak\hfil\penalty50{\tt#1}}}

\def\ti{-\allowhyphens}

\newcommand{\thismonth}{\ifcase\month 
  \or January\or February\or March\or April\or May\or June%
  \or July\or August\or September\or October\or November%
  \or December\fi}
\newcommand{\thismonthyear}{{\thismonth} {\number\year}}
\newcommand{\thisdaymonthyear}{{\number\day} {\thismonth} {\number\year}}

\def\h{{\rm H}}
\def\hb{{\rm H}_{\rm b}}
\def\ehb{{\rm EH}_{\rm b}}
\def\ha{{\rm H}_{(G,K)}}
\def\hc{{\rm H}_{\rm c}}
\def\hbc{{\rm H}_{\rm cb}}
\def\ehbc{{\rm EH}_{\rm cb}}
\def\linfty{L^\infty_{\rm w*}}
\def\linftya{L^\infty_{\mathrm{w*,alt}}}
\def\la{L^\infty_{\mathrm{alt}}}
\def\cb{{\rm C}_{\rm b}}
\def\coo{\mathrm{C}_{00}}
\def\cont{\mathrm{C}}
\def\isom{\mathrm{Isom}}
\def\aut{\mathrm{Aut}}
\def\lin{\mathrm{Lin}}
\def\supp{\mathrm{supp}}
\def\EE{\mathbf{E}}
\def\FF{\mathbf{F}}
\def\HH{\mathbf{H}}
\def\NN{\mathbf{N}}
\def\OO{\mathbf{O}}
\def\PP{\mathbf{P}}
\def\RR{\mathbf{R}}
\def\SS{\mathbf{S}}
\def\ZZ{\mathbf{Z}}
\def\C{\mathscr{C}}
\def\H{\mathscr{H}}
\def\L{\mathscr{L}}
\def\P{\mathscr{P}}
\def\Sy{\mathscr{S}}
\def\T{\mathscr{T}}
\def\Z{\mathscr{Z}}
\def\geod{\mathfrak{G}}
\def\rays{\mathfrak{R}}
\def\proba{\mathcal P}
\def\p{\partial}
\def\dim{\mathrm{dim}\,}
\def\one{\mathbf{1\kern-1.6mm 1}}
\def\bu{\bullet}
\def\weak{weak-* }
\def\cat#1{{CAT(#1)}}
\def\id{{\it I\! d}}
\def\ro{\varrho}
\def\fhi{\varphi}
\def\teta{\vartheta}
\def\epsi{\varepsilon}
\def\ti{-\allowhyphens}
\def\lra{\longrightarrow}
\def\se{\subseteq}
\def\pe#1{{^\perp\!#1}}
\def\rep{\mathrm{Rep}}
\def\prep{\mathrm{Rep}_\mathrm{par}}
\def\chom{\mathrm{Hom}_\mathrm{c}}
\def\proj{\mathrm{Pr}}
\def\amal#1{\mathop{*}\limits_{#1}}
\def\No{N\raise4pt\hbox{\tiny o}\kern+.2em}
\def\no{n\raise4pt\hbox{\tiny o}\kern+.2em}
\def\bsl{\backslash}
\def\beq{\begin{equation}}
\def\eeq{\end{equation}}
\def\bpm{\begin{pmatrix}}
\def\epm{\end{pmatrix}}
%


\title[Trees and Hyperbolic Spaces]{Equivariant Embeddings of Trees\\ into Hyperbolic Spaces}
\author{Marc Burger}
\address{M.B.: FIM, ETHZ, Z\"urich, CH-8092, Switzerland}
\email{burger@math.ethz.ch}
\author{Alessandra Iozzi}
\address{A.I.: Department Mathematik, ETHZ, Z\"urich, CH-8092, Switzerland}
\email{iozzi@math.ethz.ch}
\author{Nicolas Monod}
\address{N.M.: University of Chicago, Chicago, IL 60637, U.\,S.\,A.}
\email{monod@math.uchicago.edu}
\thanks{N.M. partially supported by FNS grant~8220-067641 and NSF grant DMS~0204601.}

\maketitle

\section{Introduction}
\label{sec_intro}%

For every cardinal $\alpha\geq 2$ there are three complete constant
curvature model manifolds of Hilbert dimension $\alpha$: the sphere
$\SS^\alpha$, the Euclidean space $\EE^\alpha$ and the hyperbolic
space $\HH^\alpha$.  Studying isometric actions on these
spaces corresponds in the first case to studying orthogonal
representations and in the second case to studying cohomology in
degree one with orthogonal representations as coefficients.  In this
paper we address the third case and, in particular, we study isometric
actions of automorphisms groups of trees on $\HH^\alpha$.

The goal of this paper is twofold: first we exhibit, for every tree
$\T$, a one-parameter family of equivariant embeddings (with respect
to appropriate representations of $\aut(\T)$) into an infinite
dimensional hyperbolic space which are, up to rescaling, asymptotically
isometric and have convex cobounded image.  Secondly, in the case in
which the tree is regular of finite valence at least 3 and
$G<\aut(\T)$ is a closed subgroup satisfying appropriate transitivity
properties, we show that the representations constructed
above give the unique irreducible nonelementary actions of $G$ by
isometries on a hyperbolic space of appropriate infinite dimension.

Quadratic forms of finite index, and in particular of index one, can
be studied on real vector spaces of arbitrary dimension.  A quadratic
form of index one leads, \emph{via} its cone of negative vectors, to a
geodesic \cat{-1} space which is then complete if and only if the
quadratic form satisfies a strong nondegeneracy condition.  
For any dimension $\alpha$, there is one such space $\HH^\alpha$ with ideal
boundary $\p\HH^\alpha$ and bordification
$\overline\HH^\alpha=\HH^\alpha\cup\p\HH^\alpha$.  Then we have:

\begin{thmintro}
\label{thm_constr_intro}%
Let $V$ be the set of vertices of a tree $\T$ with $|V|=\alpha+1$.
Then for every $\lambda>1$ there is an embedding
$\Psi_\lambda:V\to\HH^\alpha$ and a representation
$\pi_\lambda:\aut(\T)\to\isom(\HH^\alpha)$ such that:
\begin{enumerate}
\item[(i)] The map $\Psi_\lambda$ is $\pi_\lambda$\ti equivariant 
  and extends equivariantly to a boundary map
  $\p\Psi_\lambda: \p \T\to\p\HH^\alpha$ which is a homeomorphism onto
  its image.
\item[(ii)] For any two vertices $x,y\in V$ there is a precise relation
  between the combinatorial distance $d_\T$ and the Riemannian
  distance $d_{\HH^\alpha}$:
\begin{equation*}
\lambda^{d_\T(x,y)} = \cosh d_{\HH^\alpha}(\Psi_\lambda x, \Psi_\lambda y)\,.
\end{equation*}
\item[(iii)] The set $\Psi_\lambda(V)$ has finite codiameter in the convex
  hull $\C\se \HH^\alpha$ of the image of $\p\Psi_\lambda$.
\end{enumerate}\end{thmintro}

This result is the outcome of our attempt to understand 
certain claims of Gromov~\cite[Section~6.A]{Gromov91}, to the extent that 
nontrivial amalgams admit actions with unbounded orbits 
on infinite dimensional hyperbolic spaces.

Theorem~\ref{thm_constr_intro} applies in particular to 
the automorphism group $\aut(\T_r)$ of an $r$\ti regular tree $\T_r$.
By taking products and by denoting by $\HH^\infty$ the real hyperbolic
space of countable dimension $\alpha = \aleph_0$, we obtain a family of metrically
proper convex cobounded actions of $\aut(\T_r\times\T_s)$ on $\HH^\infty\times\HH^\infty$.
This leads immediately to the following:

\begin{corintro} Any cocompact lattice 
$\Gamma<\aut(\T_r\times\T_s)$ admits a metrically proper
convex cobounded action on the product $\HH^\infty\times\HH^\infty$ 
of two hyperbolic spaces of countable dimension.
\end{corintro}

Recall that in~\cite{Burger_Mozes_IHES_1} and~\cite{Burger_Mozes_IHES_2} these
type of lattices were studied systematically and examples 
of torsion-free simple groups $\Gamma$ were obtained.
These $\Gamma$'s are then fundamental groups of finite aspherical
(two-dimensional) complexes;  on the other hand,
the question to which extent there are compact aspherical manifolds
(with or without boundary) with simple fundamental group is open.
Here we obtain metrically proper convex cobounded actions of $\Gamma$
on $\HH^\infty\times\HH^\infty$;  in particular,
$\Gamma\backslash(\HH^\infty\times\HH^\infty)$ retracts to a convex
bounded (infinite dimensional) aspherical manifold with boundary.
In contrast with the algebraic aspect of this situation,
observe that if $\Lambda$ is a group acting in a metrically proper
convex cobounded way on $\HH^\infty$, 
then $\Lambda$ is a nonelementary Gromov hyperbolic group,
and hence ${\rm SQ}$\ti universal~\cite{Olshanskii_95}; in particular,
it admits many normal subgroups.

\bigskip

Turning to the classification of isometric actions, we recall that
an action of a group $G$ on $\HH^\alpha$ by isometries is {\it
  elementary} if it preserves a point in $\overline\HH^\alpha$ or a
geodesic.  Hence the study of elementary actions on $\HH^\alpha$
reduces essentially to the the study of isometric actions on the
Euclidean space $\EE^{\alpha-1}$ or on the sphere $\SS^{\alpha-1}$.
(Conversely, as observed by Gromov~\cite[7.A]{Gromov91}, any isometric action on 
$\EE^{\alpha-1}$ can be extended to an \emph{elementary} action
on $\HH^\alpha$ by realizing $\EE^{\alpha-1}$ as a horosphere based at
a fixed point at infinity; similarly, isometric actions on $\SS^{\alpha-1}$ 
extend obviously to actions with a fixed point in $\HH^\alpha$.)

Since any nonelementary action admits a unique minimal $G$\ti
invariant hyperbolic subspace (Proposition~\ref{prop_actions_alt}), we
say that a (nonelementary) action $G\to\isom(\HH^\alpha)$ is
\emph{irreducible} if there is no $G$\ti invariant hyperbolic subspace
other than $\HH^\alpha$.

\smallskip

Let now $\T_r$ be the regular tree of finite valence $r\geq3$,
$G<\aut(\T_r)$ a closed subgroup and $\pi:G\to\isom(\HH^\alpha)$ a
nonelementary action of $G$ on $\HH^\alpha$. If $G$ acts triply
transitively on $\p\T_r$, we shall see that the image $\pi(g)$ of any
hyperbolic automorphism $g\in G$ is a hyperbolic isometry of
$\HH^\alpha$ of translation length $\ell_\pi$ independent of $g$,
provided $g$ is of translation length one in $\T_r$. With this in
mind, we can state the following:

\begin{thmintro}
\label{thm_class_intro}%
Let $G<\aut(\T_r)$ be a closed subgroup which acts triply transitively
on $\p\T_r$.  For every $\ell>0$ there exists, up to equivalence, a
unique irreducible nonelementary continuous homomorphism
$\pi:G\to\isom(\HH^{\infty})$ with $\ell_\pi=\ell$.
\end{thmintro}

It follows that the representation $\pi$ in
Theorem~\ref{thm_class_intro} is exactly the irreducible component of
$\pi_\lambda|_G$ for $\lambda=e^{\ell_\pi}$ in
Theorem~\ref{thm_constr_intro}.

The structure and unitary representation theory of closed subgroups
of $\aut(\T_r)$ with some transitivity conditions
on their action at infinity is the object of intensive study.
We refer to~\cite{Figa-Talamanca_Nebbia},~\cite{Amann}, and the references therein
for a more comprehensive picture.  A first notable set of examples
to which Theorem~\ref{thm_class_intro} applies is given by the topological
group $G=\mathbf{PGL}_2(k)$, where $k$ is a non-Archimedean local
field; indeed, if $q$ is the cardinality of the residue field 
of $k$, then the action of $\mathbf{PGL}_2(k)$ on the associated 
Bruhat--Tits tree $\T_{q+1}$ identifies it with a closed subgroup 
of $\aut(\T_{q+1})$ which acts triply transitively on $\partial\T_{q+1}$.
  
Another important class of examples of closed subgroups of $\aut(\T_r)$
are the universal groups introduced in~\cite{Burger_Mozes_IHES_1}.
Recall that when $\T_r=(X,Y)$ is a $r$-regular tree, one can
label its edges in such a way that for every vertex the edges 
issued from it are labelled $\{1,2,\dots,r\}$.  Thus, for any $g\in\aut(\T_r)$
and vertex $x\in X$, one obtains a permutation $c(g,x)\in S_r$
representing $g$ ``locally'' at $x$.  To a permutation group $F<S_r$
one can then associate $U(F)$, the closed subgroup of $\aut(\T_r)$
consisting of all $g\in\aut(\T_r)$, such that $c(g,x)\in F$ for all $x\in X$.
Then $U(F)$ does not depend, up to conjugation, on the labelling of the edges.
It acts transitively on $X$ and at every vertex it induces 
the full permutation group $F$ on the edges issued from $x$ and
is, by construction, maximal with respect to this property.  
The group $U(F)$ satisfies Tits' independence condition~\cite{Tits}
and, in fact, all closed vertex transitive subgroup of $\aut(T_r)$ satisfying 
Tits' independence condition are of the form $U(F)$.  

Many properties of $U(F)$ can be read off the finite permutation group $F<S_r$.
For example, for $n=2$ and $3$, $U(F)$ is $n$-transitive on $\partial\T_r$ 
if and only if $F$ is $n$-transitive.  In the case in which $F$ is  doubly transitive, 
the unitary dual of $U(F)$ has been determined by O.~Amann~\cite{Amann}.  
When $F$ is  triply transitive, the above Theorem~\ref{thm_class_intro}
applies to $U(F)$.

\begin{rem*}
We have been informed by A.~Valette that the algebraic part of the construction of Theorem~\ref{thm_constr_intro} can also be derived elegantly from the ``tree cocycles'' that he proposes in~\cite{Valette}.
\end{rem*}

The structure of the paper is as follows.
In Section~\ref{sec_quadratic}, modelling on the finite dimensional
case, we discuss basic properties of quadratic forms of finite index
on a real vector space of arbitrary dimension, we single out the
notion of strongly nondegenerate form and show that strongly
nondegenerate forms are determined by their signature. In
Section~\ref{sec_hyperbolic} we associate a hyperbolic space to every
nondegenerate quadratic form of index one; this is a geodesic
\cat{-1} space which is complete if and only if the form is strongly
nondegenerate. This leads to the existence and the uniqueness 
of $\HH^\alpha$ for every cardinal $\alpha$.  
In Section~\ref{sec_el_actions} we discuss the
existence of irreducible hyperbolic subspaces
(Proposition~\ref{prop_actions_alt}) and establish the description of
elementary actions in terms of orthogonal representations and cocycles
in degree one (Proposition~\ref{prop_reps_isom}).  In
Section~\ref{sec_HNN},~\ref{sec_2_trans_par} and~\ref{sec_rep} we turn
more specifically to the study of actions on $\HH^\alpha$ of
automorphisms groups of trees. 
In Sections~\ref{sec_HNN} and~\ref{sec_2_trans_par} we study
more closely actions on $\HH^\alpha$ of certain locally compact
groups occurring as stabilizers of ends of trees, 
\emph{i.e.} topological ascending HNN-extensions.
These actions turn out to be elementary and hence substantial use of
Section~\ref{sec_el_actions} is made.
In Section~\ref{sec_rep} the proof of
Theorem~\ref{thm_class_intro} is completed by showing that 
the irreducible part of the action under consideration is 
determined by its restriction to
any parabolic subgroup.  In Section~\ref{sec_constr} we give the
explicit construction in Theorem~\ref{thm_constr_intro}.  Finally, the
Appendix contains some of the explicit matrix representations used
throughout the paper.

\section{Quadratic Forms of Finite Index}
\label{sec_quadratic}%

A \emph{quadratic space} is a pair $(\H,Q)$ consisting of a real
vector space $\H$ and a quadratic form $Q:\H\to \RR$. As usual, $Q$ is
\emph{positive definite} if $Q(x)>0$ for all $x\neq 0$ and \emph{negative
definite} if $-Q$ is positive definite; $\dim\H$ denotes the cardinal
of any $\RR$\ti basis of $\H$. Define
$$i_\pm(Q) = \sup\Big\{\dim W : W \text{ is a subspace of } \H \text{ and }
Q|_W \text{ is pos./neg. definite }\Big\}$$
and the {\it index} of $Q$ as 
$$i(Q) = \sup \Big\{\dim W : W \text{ is an isotropic subspace of } \H\Big\}.$$
Let $B:\H\times \H\to\RR$ be the bilinear (symmetric) form associated
to $Q$. For any subset $S\se \H$ write
$$\pe S = \big\{x\in\H : B(x,s)=0\ \forall\,s\in S\big\}.$$
We say that $Q$ is \emph{nondegenerate} if $\pe\H =0$ and that the
quadratic space is of \emph{finite index} if $i(Q)\in\NN$ (we agree
that $0\in\NN$).

Just like in the case of finite dimensional quadratic spaces, we have:

\begin{prop}
\label{prop_quad}%
Let $(\H,Q)$ be a nondegenerate quadratic space of finite index. Then
\begin{enumerate}
\item[(i)] $i(Q) = \min\{i_-(Q), i_+(Q)\}$.
\end{enumerate}
Assume now $i(Q) = i_-(Q).$
\begin{enumerate}
\item[(ii)] If $W_-\se\H$ is a negative definite subspace with
  $\dim{W_-} = i(Q)$, then $W_+:= \pe{W_-}$ is positive definite and
  $\H = W_- \oplus W_+$.
\item[(iii)] If $\H = W'_- \oplus W'_+$ is an orthogonal direct sum with
  $W'_\pm$ pos./neg. definite, then $\dim{W'_-} = i(Q)$.
\end{enumerate}
\end{prop}

We precede the proof of the proposition by a couple of lemmas.

\begin{lemma}
\label{lemma_quad1}%
Let $(\H, Q)$ be a quadratic space and $W\se \H$ a finite dimensional
subspace such that $Q|_W$ is nondegenerate. Then $\H = W\oplus \pe
W$. If moreover $Q$ is nondegenerate then $Q|_{\pe W}$ is so too.
\end{lemma}

\begin{proof}
  The kernel of the linear map $\H\to W^*$ induced by $B$ is $\pe W$;
  since $W$ has finite dimension we have
\begin{equation}
\label{eq_2_1}
\dim{(\H/\pe W)} \leq \dim{W^*} = \dim W\,.
\end{equation}
On the other hand, since $W\cap \pe W = 0$, the canonical projection
$W\to\H/\pe W$ is injective; hence it is an isomorphism
by~(\ref{eq_2_1}).
\end{proof}

\begin{lemma}
\label{lemma_quad2}%
Let $(\H, Q)$ be a quadratic space with $Q(x)\geq 0$ for all $x\in\H$.
Then $\pe\H=\{x\in \H : Q(x)=0\}$.
\end{lemma}

\begin{proof}
If $Q(x)=0$, then for all $y\in\H$ and all $\lambda\in\RR$ we have
$$0\leq Q(\lambda x+y) = 2\lambda B(x,y) + Q(y)$$
hence $B(x,y)=0$ for all $y$.
\end{proof}

\begin{proof}[Proof of Proposition~\ref{prop_quad}]
  Let $A_\pm\se \H$ be pos./neg. definite subspaces of finite
  dimension and set $A=A_- + A_+$. Then $i(Q|_A)\leq i(Q)$ and the
  theory of finite dimensional quadratic spaces implies
$$\min\{\dim{A_-}, \dim{A_+}\} \leq i(Q|_A) \leq i(Q),$$
whence $\min\{i_-(Q), i_+(Q)\} \leq i(Q)$. Assume without loss of
generality that $i_-(Q)\leq i_+(Q)$, pick a negative definite subspace
$W_-$ of dimension $i_-(Q)$ and let $W_+:= \pe{W_-}$. Since $Q|_{W_-}$
is nondegenerate, $\H = W_-\oplus W_+$ and $Q|_{W_+}$ is
nondegenerate (Lemma~\ref{lemma_quad1}). Since $\dim{W_-} = i_-(Q)$,
we have $Q(x)\geq 0$ for all $x\in W_+$ and hence, by
Lemma~\ref{lemma_quad2}, $W_+$ is positive definite. If now $W$ is an
isotropic subspace with $\dim W = i(Q)$, then $W\cap W_+ = 0$ and thus
the canonical projection $W\to \H/W_+ \cong W_-$ is injective; hence
$i(Q)\leq i_-(Q)$. This proves~(1) and~(2). As for~(3), if $W_-$ is
negative definite with $\dim{W_-} = i(Q)$, then $W_-\to \H/W'_+ \cong
W'_-$ is injective and hence an isomorphism since $\dim{W'_-} \leq
i_-(Q) = i(Q)$.
\end{proof}

In view of Proposition~\ref{prop_quad} 
we call \emph{$\pm$\ti decomposition} of $(\H,Q)$ any orthogonal
direct sum decomposition $\H = W_- \oplus W_+$ where $W_\pm$ are
pos./neg. definite. We associate to such a decomposition the scalar
product $\langle\,\, ,\,\rangle_\pm$ defined for $x,y\in\H$ by
$$\langle x,y\rangle_\pm:= B(x_+, y_+) - B(x_-, y_-)$$
where $x=x_- + x_+, y=y_- + y_+$ are the corresponding decompositions.

\begin{lemma}
\label{lemma_scalar_pm}%
Let $(\H,Q)$ be a nondegenerate quadratic space of finite index and
$\H= W_- \oplus W_+ = W'_- \oplus W'_+$ two $\pm$\ti decompositions.
Then $(\H, \langle\,\, ,\,\rangle_\pm)$ is a Hilbert space
if and only if $(\H, \langle\,\, ,\,\rangle'_\pm)$ is a Hilbert space,
in which case the two scalar products are equivalent.
\end{lemma}

We need the following:

\begin{lemma}
\label{lemma_scalar_codim}%
Let $\H$ be a real vector space, $\langle\,\, ,\,\rangle_1,\, \langle
\,\,,\,\rangle_2$ two scalar products and $\H'\se \H$ a subspace such
that
\begin{enumerate}
\item $\langle \,\,,\,\rangle_1,\, \langle\,\, ,\,\rangle_2$ coincide
  on $\H'$ and $\H'$ is complete;
\item $\H'$ is of finite codimension.
\end{enumerate}
Then $\langle \,\,,\,\rangle_1,\, \langle\,\, ,\,\rangle_2$ are
equivalent and $\H$ is a Hilbert space.
\end{lemma}

\begin{proof}
  Let $\H'_1$ be the orthogonal of $\H'$ for $\langle
  \,\,,\,\rangle_1$. Since $\H'$ is complete, we have $\H = \H' \oplus
  \H'_1$. But $\H'_1$ is complete because it is finite dimensional and
  hence $(\H, \langle \,\,,\,\rangle_1)$ is a Hilbert space. For any
  $x\in \H$, write $x = x' + x'_1$ according to the above
  decomposition. Since $\|\,\|_1$ and $\|\,\|_2$ are equivalent on
  $\H'_1$, there is $c>0$ with $\|x'_1\|_1^2 \geq c \|x'_1\|_2^2$ for
  all $x$. We may chose $c\leq 1$ and now
\begin{multline*}
  \|x\|_1^2 = \|x'\|_1^2 + \|x'_1\|_1^2 \geq \|x'\|_1^2 + c\|x'_1\|_2^2 \geq \\
  \geq c \big(\|x'\|_1^2 + \|x'_1\|_2^2\big) \geq \frac{c}{2}
  \big(\|x'\|_1 + \|x'_1\|_2\big)^2 \geq \frac{c}{2} \|x\|_2^2.
\end{multline*}
\end{proof}

\begin{proof}[Proof of Lemma~\ref{lemma_scalar_pm}]
  Assume $i(Q) = i_-(Q)$. Since $B$ is continuous with respect to both
  $\|\,\|_\pm$ and $\|\,\|'_\pm$, all subspaces considered are
  closed for both topologies and so is in particular $W_+\cap W'_+$.
  Moreover, the latter is of codimension at most $2i(Q)$, hence we
  conclude by Lemma~\ref{lemma_scalar_codim}.
\end{proof}

\begin{defi}
  A nondegenerate quadratic space of finite index $(\H,Q)$ is
  \emph{strongly nondegenerate} if for some (hence any) $\pm$\ti
  decomposition $\H=W_-\oplus W_+$ the space $(\H, \langle
  \,\,,\,\rangle_\pm)$ is a Hilbert space.
\end{defi}

We denote by $d(\L)$ the cardinal of any Hilbert basis of a Hilbert
space $\L$. Observing that $d(\L)$ depends only on the equivalence
class of the scalar product, we deduce that the pair $(d(W_+),
d(W_-))$ is independent of the choice of a $\pm$\ti decomposition $\H
= W_- \oplus W_+$. We call $(d(W_+), d(W_-))$ the \emph{signature} of
the strongly nondegenerate quadratic space $(\H, Q)$.

Two quadratic spaces $(\H_1, Q_1)$ and $(\H_2, Q_2)$ are
\emph{isomorphic} if there is a vector space isomorphism $T: \H_1 \to
\H_2$ with $Q_1 = Q_2\circ T$. Observe that if $(\H_1, Q_1)$ is
nondegenerate of finite index, then so is $(\H_2, Q_2)$. Since the
image of a $\pm$\ti decomposition for $Q_1$ is a $\pm$\ti
decomposition for $Q_2$, we see that $Q_2$ is strongly nondegenerate
if $Q_1$ is so; in that case $T$ is automatically continuous. In
particular the orthogonal group $\OO(Q)$ of a strongly nondegenerate
form $Q$ of finite index consists of bounded linear operators.

\begin{prop}
\label{prop_signature}%
For strongly nondegenerate forms of finite index, the signature is a
complete invariant of isomorphisms.
\end{prop}

\begin{proof}
  Using $\pm$\ti decompositions, this follows immediately from the
  fact that $d(\L)$ determines completely the Hilbert   spaces $\L$
  up to isomorphisms.
\end{proof}

Let $(\H, Q)$ be a strongly nondegenerate form of finite index,
$\H^*$ the topological dual and $A: \H \to \H^*$ the continuous
morphism associated to $B$. Applying the Riesz representation theorem
to the restrictions of $B$ to $W_\pm$ for a $\pm$\ti decomposition $\H
= W_- \oplus W_+$, we deduce that $A$ is an isomorphism (of
topological vector spaces).

\begin{prop}
\label{prop_nd_subspace}%
Let $(\H, Q)$ be a strongly nondegenerate form of finite index and
$V\se \H$ a closed subspace such that $Q|_V$ is nondegenerate. Then
$(V,Q|_V)$ is strongly nondegenerate and $\H = V\oplus \pe V$.
\end{prop}

\begin{proof}
  Assume $i(Q) = i_-(Q)$ and let $V = U_- \oplus U_+$ be a $\pm$\ti
  decomposition of $V$ with $U_\pm$ pos./neg. definite (which exists
  by Proposition~\ref{prop_quad}).  Since $V$ is closed and $B$
  continuous, $U_+ = \pe U_- \cap V$ is closed. Let now $\H = W_-
  \oplus W_+$ be any $\pm$\ti decomposition of $\H$. Then $W_+\cap
  U_+$ is closed, of finite codimension in $V$ and $B$ coincides with
  $\langle \,\,,\,\rangle_\pm$ on it. By
  Lemma~\ref{lemma_scalar_codim} with
  $\langle\,\,,\,\rangle_1=\langle\,\,,\,\rangle_2$, we deduce
  that $V$ is a Hilbert space and hence $(V,Q|_V)$ is strongly
  nondegenerate.  To conclude, the nondegeneracy of $(V,Q|_V)$
  implies, as observed above, that the morphism $A_V:V\to V^*$
  associated to $B|_V$ is a topological isomorphism. In particular,
  for every $x\in\H$ there is $x_V=A_V^{-1}B(x,\cdot)|_V \in V^*$ such
  that $B(x,y) = B(x_V,y)$ for all $y\in V$. Thus $x \in \pe V + x_V$
  and the claim follows.
\end{proof}

\section{Real Hyperbolic Space}
\label{sec_hyperbolic}%


Let $(\H, Q)$ be a nondegenerate quadratic space of index one; we
assume more specifically that $i(Q) = i_-(Q) = 1$. Let $C_-:= \big\{x\in \H
: Q(x) < 0 \big\}$ be the cone of negative vectors and $\HH:=\RR^*\bsl
C_-$ the set of negative lines. One proves as usual the reverse Cauchy--Schwartz inequality
%
$$B(x,y)^2 \geq Q(x)Q(y)\,,\kern1cm \forall\,x,y\in C_-$$
with equality if and only if $\RR^*x = \RR^* y$. This allows to define
$\widetilde d: C_-\times C_-\to\RR_+$ by
%
$$\cosh^2 \widetilde d(x,y)= \frac{B(x,y)^2}{Q(x) Q(y)}\,,$$ 
which descends to a well defined function $d: \HH\times \HH\to\RR_+$.

\begin{rem}\label{rem_fin_dim}
A useful geometric fact is that for any finite set $S$ of negative lines, 
the restriction $Q|_{\H_S}$ of $Q$ to the span $\H_S\se \H$ of $S$ is 
equivalent to the standard real quadratic form of signature $(\dim{\H_S}-1, 1)$ 
and that, under this isomorphism, the restriction of $d$ to the image
$\HH_S$ of $\H_S$ in $\HH$ corresponds to the standard distance on the
finite dimensional real hyperbolic space of dimension $|S|-1$.
\end{rem}
As a consequence of the above remark and~\cite[Theorem 10.10]{Bridson_Haefliger},
we have the following:

\begin{prop}%
\label{rem_fd_hyp}%
  The function $d$ is a distance function with respect to which $\HH$
  is a geodesic \cat{-1} space.\hfill\qedsymbol
\end{prop}


Let $\ell_-$ be a vector of length $-1$ and let $\H=\H_+\oplus\RR\ell_-$ 
be the orthogonal decomposition, where $Q|_{\H_+}$ is positive definite 
(see Proposition~\ref{prop_quad}).  The exponential map
$\exp: \H_+\to \HH$ is defined as follows. For every $v\in\H_+$ 
there is a unique $t>0$ such that $x:=v+t\ell_-$ has length $-1$;
we define $\exp(v):=[x]$ to be the image of $x$ in $\HH$.
A straightforward computation gives
$$\cosh d\big(\exp(v), \exp(w)\big) =\big|-B(v,w) + \sqrt{1+Q(v)}\sqrt{1+Q(w)}\big|\,,$$
and, in particular, $\cosh d(\exp(v), [\ell_-]) = \sqrt{1+Q(v)}$.

\begin{prop}
  The \cat{-1} space $\HH$ is complete if and only if $(\H, Q)$ is
  strongly nondegenerate.
\end{prop}

\begin{proof}
  Using the above formul\ae, on checks that $\exp$ is for all $R>0$ a
  bi-Lipschitz bijection between the ball in $(\H_+, Q|_{\H_+})$ of
  radius $(\sinh{R})^2$ centered at $0$ and the ball in $(\HH,d)$ of
  radius $R$ centered at $[\ell_-]$. 
  Thus $(\HH,d)$ is complete if and only if $\H_+$ is complete,
  which in view of Lemma~\ref{lemma_scalar_codim} (with
  $\langle\,\,,\,\rangle_1=\langle\,\,,\,\rangle_2$)
  is equivalent to the quadratic space $(\H,Q)$ being strongly
  nondegenerate.
\end{proof}

Observe that any orthogonal transformation $T\in\OO(Q)$ preserves
$C_-$ and descends to an isometry of $\HH$. The group $\OO(Q)$ is a
direct product $\OO_+(Q)\cdot\{\pm\id\}$, where $\OO_+(Q)$ is the
subgroup preserving the (two) connected components of $C_-$. Let
$\PP\OO(Q)=\OO(Q)/\pm\id$.  Then:

\begin{prop}
\label{prop_isom_H}%
The homomorphism $\mathbf{O}(Q)\rightarrow\isom(\HH)$ induces isomorphisms 
\begin{equation*}
\OO_+(Q)\rightarrow\mathbf{PO}(Q)\rightarrow\isom(\HH)\,.
\end{equation*}
\end{prop}

\begin{rem}
\label{rem_cont}%
Let $G$ be a topological group and $\pi: G\to \OO(Q)$ a group
homomorphism. We call $\pi$ \emph{continuous} if the action map
$G\times \H\to\H$ is continuous; then, the resulting action $G\times
\HH\to \HH$ is continuous.  Conversely, given a continuous action
$G\times\HH\to\HH$, one verifies that the resulting homomorphism
$G\to\OO_+(Q)$ deduced from Proposition~\ref{prop_isom_H} is
continuous.
\end{rem}

\begin{lemma}
\label{lemma_trans}%
The $\OO(Q)$\ti action on $\HH$ is transitive.
\end{lemma}

\begin{proof}[Proof of the lemma]
Let $L,L'\se \H$ be two negative lines. Then $L\oplus \pe{L}$ and $L'\oplus \pe{L'}$ are two $\pm$\ti decompositions, and $\pe{L}$, $\pe{L'}$ are isomorphic Hilbert spaces. Hence there is (compare also Proposition~\ref{prop_signature}) an isomorphism of $(\H,Q)$ bringing $L$ to $L'$.
\end{proof}

\begin{proof}[Proof of Proposition~\ref{prop_isom_H}]
  Let $T\in\isom(\HH)$. By Lemma~\ref{lemma_trans}, 
  we may assume that $T$ fixes $[\ell_-]$. Define a map $U:\H_+\to\H_+$ 
  by $\exp(U(v)) = T(\exp(v))$. It follows from the above formul\ae\ 
  that $U$ is a bijection, fixes $0$ and preserves $B|_{\H_+}$. Hence $U$ is a 
  linear orthogonal transformation of $\H_+$. Defining $S:=U\oplus \id$,
  one verifies that $S\in\OO_+(Q)$ corresponds to $T$ \emph{via} $\OO(Q)\to\isom(\HH)$ 
  and the statement follows.
\end{proof}

One proves similarly: 

\begin{prop}
  Let $(\H_i, Q_i)$ be strongly nondegenerate quadratic spaces of
  signature $(\alpha_i, 1)$ and let $\HH_i$ be the associated
  hyperbolic spaces (for $i=1,2)$. The following are equivalent:
\begin{enumerate}
\item[$\bullet$] $(\H_1, Q_1)$ is isomorphic to $(\H_2, Q_2)$.
\item[$\bullet$] $\alpha_1 = \alpha_2$.
\item[$\bullet$] $\HH_1$ is isometric to $\HH_2$.\hfill\qedsymbol
\end{enumerate}
\end{prop}

Thus we obtain for each cardinal $\alpha$ ``the'' real hyperbolic
space $\HH^\alpha$. 

\subsection{Bordification}

Let again $\HH = \RR^*\bsl C_-$ be the real hyperbolic space
associated to a strongly nondegenerate quadratic space $(\H, Q)$ of
signature $(\alpha, 1)$. Let $\p\HH$ be the boundary of the \cat{-1}
space $(\HH,d)$ defined as usual as classes of asymptotic rays. Set
$$C_0:= \big\{x\in\H : Q(x) =0, x\neq 0\big\}, \kern1cm C_{\leq 0} :=
\big\{x\in\H : Q(x) \leq 0, x\neq 0\big\}.$$
Using that any configuration of finitely many geodesics in $\HH$ is
contained in a finite dimensional hyperbolic subspace (see
Remark~\ref{rem_fin_dim}),
we obtain a bijection identifying the
bordification $\overline{\HH} = \HH\sqcup \p \HH$ with the set
$\RR^*\bsl C_{\leq 0}$. We relate this to the description of
$\overline{\HH}$ in terms of Busemann cocycles: for every $x\in
C_{\leq0}$, define $\widetilde b_x:C_-\times C_-\to\RR$ by
\begin{equation}\label{eq_buse}
\widetilde b_x(y,z):=
\begin{cases}
\widetilde d(x,y) - \widetilde d(x,z)&\text{ if }\,x\in C_-,\\
\frac{1}{2} \ln\frac{B(x,y)^2 Q(z)}{B(x,z)^2 Q(y)}&\text{ if }x\in C_0\,.
\end{cases}
\end{equation}
%
Then, for every $x\in C_{\leq0}$, $\widetilde b_x$ satisfies the
cocycle identity
$$\widetilde b_x(y_2, y_3) - \widetilde b_x(y_1, y_3) + \widetilde b_x(y_1, y_2) = 0.$$
Moreover, $\widetilde b_x$ gives a well defined function $b:
\RR^*\bsl C_{\leq 0} \times \HH\times \HH \to\RR$ which coincides on
$\HH\times \HH\times \HH$ with $(x,y,z)\mapsto d(x,y) - d(x,z)$.

%
%

For every $x\in\RR^*\bsl C_{\leq 0}$ and $z_1,z_2\in\HH$,
the cocycle property implies that the continuous functions
$y\mapsto b_x(y,z_1)$ and $y\mapsto b_x(y,z_2)$, defined on $\HH$,
differ by a constant.  Thus we obtain, for every $x\in\RR^*\bsl C_{\leq 0}$,
a well defined class $B(x)\in\cont(\HH)/\RR$,
where $\cont(\HH)$ is the space of continuous functions on $\HH$.
Endowing $\cont(\HH)/\RR$ with the topology coming 
from the topology on $\cont(\HH)$ of uniform convergence 
on bounded sets, and denoting by $\overline{B(\HH)}$ the closure of
$B(\HH)$ in $\cont(\HH)/\RR$, we have the following:

\begin{prop}
  The map $B:\RR^*\bsl C_{\leq 0} \to
  \overline{B(\HH)}$ is a homeomorphism when 
  $\RR^*\bsl C_{\leq 0}$ is endowed with the quotient of the norm topology.
\end{prop}

One can verify that in this topology $\overline\HH$ is compact if and
only if $\HH$ is finite dimensional.

By a slight abuse of terminology, we call {\it horospheres} in $\H$ 
(respectively, in $\HH$) centered at $x\in C_0$ (respectively, at $\xi\in\p\HH$) 
the level sets of $\tilde b_x(\cdot,z)$ (respectively, $B(\xi)$).  
%

\section{Nonelementary and Elementary Actions}
\label{sec_el_actions}%

In this section we study basic properties of group actions 
on hyperbolic spaces.  First we establish that any nonelementary 
action has a unique minimal invariant hyperbolic subspace, 
and then we turn to the description of elementary actions 
in terms of orthogonal representations, characters, and continuous
cocycles.

\medskip

Let $X$ be a metric space. Recall that a \emph{semicontraction} is a
map $T:X\to X$ such that $d(Tx, Ty)\leq d(x,y)$ for all $x,y\in X$. Recall the

\begin{prop}
\label{prop_karlsson}%
Let $X$ be a complete \cat{-1} space and $T: X\to X$ a
semicontraction. Then one of the following holds:
\begin{enumerate}
\item[(i)] The set $\{T^n x : n\geq 1\}$ is bounded for some (hence any)
  $x\in X$ and the set $X^T\se X$ of $T$-fixed points is not
  empty.
\item[(ii)] The set $\{T^n x : n\geq 1\}$ is unbounded for some (hence any)
  $x\in X$ and there exists a subsequence $\{n_k\}_{k\geq 1}$ and
  $\xi\in\p X$ with $\lim_{k\to\infty} T^{n_k} x = \xi$ and $T \xi =
  \xi$. Moreover, $|(\p X)^T| = 1\text{ or } 2$.
\end{enumerate}
\end{prop}

(A general semicontraction need not extend to infinity; the notation $T \xi = \xi$ means that $T: X\to X$ extends by continuity to $X\cup\{\xi\}$ endowed with the topology induced by $\overline{\! X}$.)

\begin{proof}[Proof of the proposition]
The case~(ii) follows from the argument given by A.~Karlsson (proof of~\cite[5.1]{Karlsson99}); see~\cite[\S3]{Karlsson-Noskov} for the additional statement on $|(\p X)^T|$. If on the other hand the $T$\ti orbits are bounded, then it is known that $X^T$ is non-empty; indeed, one verifies that for any $x\in X$ the circumcentre of the set $\{T^k(x) : k\geq n\}$ converges to a $T$\ti fixed point as $n\to \infty$.
\end{proof}

This result applies in particular to the case where $T$ is an isometry and is the basis for the classification of isometries.

\begin{defi}
An isometry $T$ is called:
\begin{enumerate}
\item[--] \emph{Elliptic} if $\{T^n x : n\geq 1\}$ is bounded.
\item[--] \emph{Parabolic} if $\{T^n x : n\geq 1\}$ is unbounded and $|(\p X)^T| = 1$.
\item[--] \emph{Hyperbolic} if $\{T^n x : n\geq 1\}$ is unbounded and $|(\p X)^T| = 2$.
\end{enumerate}
\end{defi}

If $T$ is hyperbolic, then $(\p X)^T = \{\xi_-, \xi_+\}$ with $\lim_{n\to\pm\infty}
T^n x = \xi_\pm$ for all $x\in X$ (see again~\cite[\S3]{Karlsson-Noskov}). However
, when $X$ is not proper and $T$ is
parabolic, the sequence $T^n x$ might not converge in $\overline{X}$.

\medskip

As usual a group action $G\times X\to X$ by isometries is called
\emph{elementary} if $G$ preserves a nonempty finite subset of
$\overline{X}$. This is equivalent to saying that either $G$ fixes a
point in $\overline{X}$ or it preserves a geodesic.

\bigskip

Let now $\HH$ be the hyperbolic space associated to a strongly
nondegenerate quadratic space $(\H, Q)$ of signature $(\alpha, 1)$.
In the sequel we shall study nonelementary and elementary actions 
and we shall prove the following

\begin{prop}
\label{prop_actions_alt}%
Let $\pi: G\to\OO(Q)$ be a homomorphism. Then, one of the following holds:
\begin{enumerate}
\item[(i)] $G$ preserves an isotropic line and all horospheres in $\H$
  centered at it.
\item[(ii)] $G$ preserves a negative line.
\item[(iii)] There is a unique minimal nondegenerate closed $G$\ti invariant
  subspace $\H_1\se\H$ of index one. Any nondegenerate closed $G$\ti
  invariant subspace of index one contains $\H_1$.
\end{enumerate}
\end{prop}

First we prove the proposition in the case where the associated action
on $\HH$ is nonelementary, which excludes of course (i) and (ii). The
remaining will be a consequence of a closer analysis of elementary
actions.

\begin{proof}[Proof in the nonelementary case]
  We need to show that (iii) holds. Let $\P$ be the set of $G$\ti
  invariant closed positive definite subspaces of $\H$, ordered by
  inclusion, let $\C\se\P$ be a maximal chain and $L:=
  \overline{\bigcup \C}$. Then $L$ is closed, $G$\ti invariant, and
  $Q|_L \geq 0$. By Lemma~\ref{lemma_quad2} applied to $L$, we have
  that $\{x\in L : Q(x) = 0\} = L \cap \pe L$; if the latter were not
  zero, it would be a $G$\ti invariant isotropic line, contradicting
  the assumption that the action is nonelementary. 
  Thus $L$ is a $G$\ti invariant closed maximal
  positive definite subspace of $\H$. Set $\H_0:=L$, $\H_1 := \pe L$. 
  Then, by Proposition~\ref{prop_nd_subspace}, $\H = \H_0 \oplus \H_1$ 
  and $\H_1$ is a $G$\ti invariant closed nondegenerate subspace of 
  index one which is minimal with respect to these properties.
  
  Let now $\H = \H'_0 \oplus \H'_1$ be any other orthogonal
  decomposition into $G$\ti invariant closed subspaces where $\H'_0$
  is positive definite and $\H'_1$ of index one. We need to show that
  $\H'_1\supseteq\H_1$. Consider $J:=\H_1 \cap \H'_1$. Again, 
  since the $G$\ti action is not elementary, $J$ is
  nondegenerate. There are two cases:

\smallskip

{\it $J$ is indefinite.} Since $J\se\H_1$, we have either $J=\H_1$,
whence $\H'_1\supseteq \H_1$ and we are done; or $J=0$.

{\it $J$ is positive definite.} Then $\H_0\oplus J$ would be a $G$\ti
invariant closed positive definite subspace and, by maximality, we
would have $\H_0 \oplus J\se \H_0$ and hence, once again, $J=0$.

Thus, we may assume (for a contradiction) $\H_1\cap \H'_1 = 0$. Let
$\HH_1, \HH'_1\se\HH$ be the corresponding hyperbolic subspaces and
consider the orthogonal projections $p:\HH\to\HH_1$ and
$p':\HH\to\HH'_1$ given by the nearest point retraction.  
Since $\HH_1\cap\HH'_1=\varnothing$ and $\HH$ is
\cat{-1}, both $p|_{\HH'_1}$ and $p'|_{\HH_1}$ are contractions.  Hence
the map $f: \HH_1\to \HH_1$, defined by $f:= p|_{\HH'_1}\circ
p'|_{\HH_1}$, is a $G$\ti equivariant contraction, that is, $d(f(x),
f(y)) < d(x,y)$ for all distinct $x,y\in \HH_1$. If for some
$x\in\HH_1$ the set $\{d(f^n(x),x):\,n\geq1\}$ were bounded,
Proposition~\ref{prop_karlsson}(i) would imply the existence of an
$f$\ti fixed point in $\HH$.  Since $f$ is $G$\ti equivariant, the set
of its fixed points $\HH_1^f$ is $G$\ti invariant.  However,
since $f$ is a contraction, $\HH_1^f$ consists of one point, which is
hence $G$\ti fixed, contradicting the assumption that the $G$\ti action
is nonelementary.

Thus, the $f$\ti orbits are unbounded and by
Proposition~\ref{prop_karlsson}(ii) there is a subsequence $\{n_k\}$
and $\xi\in\p\HH_1$ with $\lim_{k\to\infty} f^{n_k} (x) = \xi$. However,
for every given $g\in G$ and for $x\in\HH_1$
$$d(g f^n(x), f^n(x)) = d(f^n(g x), f^n(x)) < d(gx, x)$$
is bounded independently of $n$, thus, by passing to subsequences, the
sequences $gf^n(x)$ and $f^n(x)$ are at bounded distance and hence
define the same point at infinity, namely $g\xi=\xi$.  Since this 
contradicts again the assumption that the action is nonelementary, 
the proof in this case is complete.
\end{proof}

\subsection{Elementary Actions}
\label{sec_elem}%

As before, let $(\H, Q)$ be a strongly nondegenerate quadratic space of signature 
$(\alpha, 1)$. We shall study, for a topological group $G$, 
the elementary actions on $\HH^\alpha$, and more specifically, 
the actions fixing a point in $\partial\HH^\alpha$.  
Thus, fix $L_+$ a isotropic line in $\H$, let $\OO_{L_+}(Q)$
be its stabilizer in $\OO(Q)$, and let $\rep(G, \OO_{L_+}(Q))$
be the set of continuous representations (see Remark~\ref{rem_cont}).
Then $G$ acts on $L_+$ by multiplication by a continuous
character $\chi:G\to \RR^*$. The bilinear form $B$ induces a Hilbert
space structure of Hilbert dimension $\alpha-1$ on $\pe{L_+}/L_+$ and
we may thus fix a real Hilbert space $E$ of Hilbert dimension
$\alpha-1$ and an isomorphism $i:\pe{L_+}/L_+\to E$. Since
$\pi:G\to\OO_{L_+}$ is continuous and $G$ preserves $L_+$ and
$\pe{L_+}$, it induces on $\pe{L_+}/L_+$ a continuous
orthogonal representation which we transport to $E$ \emph{via} $i$,
obtaining $\ro: G\to \OO(E)$. The space $\H/\pe{L_+}$ is
one-dimensional; since $G$ preserves $B$, it acts on that space by
multiplication by $\chi^{-1}$.

We define a new $G$\ti module structure on $\H$ by means of the
continuous representation $\chi\otimes \pi$. Then we have a short
exact sequence of $G$\ti modules
$$0\lra \pe{L_+}/L_+\lra \H/L_+ \lra \H/\pe{L_+}\lra 0$$
in which the last term is a trivial $G$\ti module of dimension one.
Thus, applying the corresponding transgression map in degree zero, the
image of the trivial module $\H/\pe{L_+}$ in the continuous cohomology
$\hc^1(G,\pe{L_+}/L_+)$ yields \emph{via} $i$ a subspace
$$\RR \cdot\eta \se \hc^1(G,\chi\otimes\ro)\,.$$
Thus, with $i$ fixed, we associated 
to every continuous homomorphism $\pi:G\to\OO_{L_+}(Q)$
the following data:
\begin{itemize}
\item[--] a continuous homomorphism $\chi_\pi\in\chom(G,\RR^*)$,
\item[--] a continuous orthogonal representation $\ro_\pi\in\chom(G, \OO(E))$, 
          and
\item[--] a continuous class $\eta_\pi\in\hc^1(G,\chi\otimes\ro)$, well
defined up to scalar multiplication.
\end{itemize}

Denoting by $\Z(G,E)$ the set of
all triples $(\chi, \ro, \eta)$ with $\chi\in\chom(G,\RR^*)$,
$\ro\in\chom(G,\OO(E))$ and $\eta\in\hc^1(G,\chi\otimes \ro)$, we have
that $\OO(E)\times \RR^*$ acts on $\Z(G,E)$ by $(T,\lambda)(\chi,
\ro,\eta) = (\chi, T\ro T^{-1}, \lambda T\eta)$.

\begin{prop}
\label{prop_reps_isom}%
\begin{enumerate}
\item[(i)] The map $\rep(G, \OO_{L_+}(Q)) \to \Z(G,E)$, $\pi\mapsto
  (\chi_\pi, \ro_\pi, \eta_\pi)$ induces a bijection
\begin{equation*}
\rep(G, \OO_{L_+}(Q)) \big/ \OO_{L_+}(Q) \xrightarrow{\ \cong\ } \big[\OO(E)\times \RR^*\big] \big\bsl\Z(G,E)\,.
\end{equation*}
\item[(ii)] The representation $\pi$ leaves all horospheres centered at
  $L_+$ invariant if and only if $|\chi_\pi| = 1$.
\end{enumerate}
\end{prop}

\begin{proof}
  Given $\chi:G\to \RR^*$, $\ro:G\to \OO(E)$ and
  $\eta\in\hc^1(G, \chi\otimes \ro)$, 
  we indicate how to reconstruct $\pi\in\rep(G,
  \OO_{L_+}(Q))$.  Fix an isotropic line $L_-\neq L_+$ and let $f:G\to
  E$ be a continuous cocycle representing $\eta$. Set $F =
  \pe{(L_+\oplus L_-)}$, and denote by $j$ the isomorphism of Hilbert
  spaces obtained by composing $F\to \pe{L_+}\to E$. Fix $\ell_\pm\in
  L_\pm$ with $B(\ell_-, \ell_+) = 1$. Using the notation of
  Appendix~\ref{sec_app}, we define
$$\pi(g) = \bpm
\pi(g)_1 & \begin{matrix} \pi(g)^+_2\\0\\ \end{matrix}\\
\begin{matrix} 0 & \pi(g)^-_3\end{matrix} & \pi(g)_4\\
\epm$$
by
$$\pi(g)_1 = \bpm
\chi(g) & a(g) \\
0 & \chi(g)^{-1}\\
\epm$$
where
\begin{equation}
\begin{aligned}
\label{eq_rep}
a(g) &= -\frac{1}{2} \chi(g) \|f(g)\|^2_E\\
\pi(g)^+_2(v) &= -\langle \ro(g) j(v), f(g) \rangle_E,\kern1cm \forall\,v\in F\\
\pi(g)^-_3 &= \chi(g)^{-1} j^{-1}(f(g))\\
\pi(g)_4 &= j^{-1}\ro(g) j.
\end{aligned}
\end{equation}
The rest of the proposition is now a verification left to the reader
and uses the fact that $|\chi_\pi|$ is the exponential of the Busemann
character associated to the fixed point $L_+$.
\end{proof}
\bigskip

We turn now to representations $\pi: G\to \OO_{L_+}(Q)$ for which
$|\chi|$ is not identically $1$ (write $\chi=\chi_\pi$). We fix once
and for all $\chi$ and $a\in G$ such that $|\chi(a)|\neq 1$.

\begin{defi}
  We shall say that a continuous cocycle $f:G\to E$ for
  $\chi\otimes\ro$ is \emph{standard} if $f(a)=0$.
\end{defi}

\begin{lemma}
\label{lemma_std_cocycle}%
Let $\ro:G\to\OO(E)$ be a continuous orthogonal
representation and $\chi:G\to\RR^\times$ a continuous homomorphism
with $|\chi(a)|\neq1$. 
\begin{enumerate}
\item[(i)] Every class in $\hc^1(G,\chi\otimes\ro)$ admits a unique
  standard representative.
\item[(ii)] If $M<G$ is a compact subgroup normalised by $a$, then any
  standard cocycle vanishes on $M$.
\end{enumerate}
\end{lemma}

\begin{proof}
  {\it (i)}~Set $\tau = \chi\otimes \ro$. Recall that $|\chi(a)|\neq 1$ 
implies that $1-\tau(a)$ is invertible and hence $\h^1(\langle a\rangle, \tau)$ 
vanishes. Therefore, for any
  cocycle $f':G\to E$ there exists $v\in E$ such that $f'(a^n) =
  \tau(a^n) v -v$. Now $f(g):= f'(g) + v- \tau(g)v$ defines a standard
  cocycle.  If $f_1$ and $f_2$ are any two cohomologous cocycles,
  there exists $v\in E$ such that $f_1(g)=f_2(g)+\tau(g)v-v$.  If in
  addition $f_1$ and $f_2$ are standard, then $\tau(a)v=v$, which
  implies, since $|\chi(a)|\neq1$, that $v=0$ and hence $f_1=f_2$.
  
  {\it (ii)}~Let $f$ be a standard cocycle. Since $M$ is compact, $C:=\sup_{k\in M}\|f(g)\|$ is finite. We
  have for all $k\in M$
$$f(k) = \tau(k) f(a) + f(k) = f(ka) = f(a a^{-1} k a) = \tau(a) f(a^{-1} k a) + f(a) = \tau(a) f(a^{-1} k a),$$
which implies that $C = |\chi(a)| C$ and hence $C=0$.
\end{proof}

Let $\eta\in\hc^1(G,\chi\otimes\ro)$ and let $\pi:G\to\OO_{L_+}(Q)$ be
the homomorphism associated to $(\chi,\ro,\eta)$ by the above
construction.  Then $\pi(a)$ is hyperbolic with fixed points $L_-$ and
$L_+$.  If $f:G\to E$ is the standard cocycle representing the class
$\eta$, define
$$I_\eta:= \overline{\langle f(g) : g\in G\rangle}\,,$$
which is a closed $G$\ti invariant subspace of $E$. 
With these definitions, we have the following:

\begin{prop}
\label{prop_dec_chi_neq_1}%
There is a $G$\ti invariant orthogonal decomposition
$$\H = \H_1 \oplus \H_0,$$
where
$$\H_1 = (L_+ \oplus L_-) \oplus j^{-1}(I_\eta),\kern1cm \H_0 = \pe{\H_1}.$$
Moreover, the subspace $\H_1$ is nondegenerate of index one, $\H_0$
is positive definite and if $\H'_1\se \H$ is any closed
nondegenerate $G$\ti invariant subspace of index one, then 
$\H'_1 \supset \H_1$.
\end{prop}

\begin{proof}
  We verify the last assertion: let $\HH'_1\se \HH_1$ be the
  hyperbolic subspace associated to $\H'_1$. Without loss of
  generality, $|\chi(a)|>1$. Then $\lim_{n\to\pm\infty} \pi(a^n) x =
  L_{\pm}$ for all $x\in\HH'_1$, hence $L_\pm\in\p \HH'_1$. If
  $L:=L_+\oplus L_-$ and $\proj_L:\H\to L$ is the orthogonal
  projection, this implies that $\pi(g) \ell_- -
  \proj_L(\pi(g)\ell_-)$ is in $\H'_1$ for all $g\in G$; therefore,
  using the formul\ae\ in (\ref{eq_rep}) and Appendix~\ref{sec_app},
  $\H'_1 \supseteq j^{-1}(I_\eta)$.
\end{proof}

\begin{proof}[End of proof of Proposition~\ref{prop_actions_alt}]
  We know already that the alternative~(iii) of the proposition holds in
  the nonelementary case. If the $G$\ti action is elementary, then
  there are the following possibilities:

{\it $G$ fixes a point in $\HH$}.  This corresponds to~(ii).

{\it $G$ leaves a geodesic line in $\HH$ invariant}. Denote this line
by $\HH_1\se\HH$. This means that there is a $G$\ti invariant
two-dimensional subspace $\H_1\subseteq\H$ of index one. In fact,
$\H_1 = L_+\oplus L_-$, where $\p \HH_1 = \{L_-, L_+\}$. We may assume
that $G$ has no fixed point in $\HH_1$, so there is $g\in G$ acting
hyperbolically: say, $\lim_{n\to\pm\infty} \pi(g)^n x = L_\pm$. But
then, if $\H'_1$ is any $G$\ti invariant closed nondegenerate
subspace of index one and $\HH'_1\se\HH$ the associated hyperbolic
subspace, we get $\p\HH'_1\ni L_\pm$, hence $\H'_1 \supseteq \H_1$.

{\it $G$ fixes an isotropic line $L_+$}. Let then $\chi:G\to \RR^*$ be
the associated character. Either $|\chi|=1$ and we are in
alternative~(i), or $|\chi|\neq 1$ and we can apply
Proposition~\ref{prop_dec_chi_neq_1}, so that alternative (iii) holds
once again.
\end{proof}

\section{Actions of Certain HNN-Extensions}
\label{sec_HNN}%

Let $P = \langle a\rangle\ltimes N$ be a locally compact group,
semidirect product of an infinite cyclic subgroup $\langle a\rangle$
with generator $a$ and a closed normal subgroup $N = \bigcup_{n\in\ZZ}
K_n$ which is the increasing union of compact open subgroups $K_n <
K_{n+1}$ such that $aK_n a^{-1} = K_{n+1}$ for all $n\in \ZZ$. 
In other words, $P$ is the topological ascending HNN-extension $P=K_0 *_a$.

In this section we study more in detail elementary actions of $P$ on
a hyperbolic space $\HH$, using the results in Section~\ref{sec_el_actions},
in particular Proposition~\ref{prop_reps_isom}.  We begin with
the following general fact:

\begin{lemma}
\label{lem_class_isom}%
Let $X$ be a complete \cat{-1} space and $P\times X\to X$ a continuous
action by isometries. Then there is a $P$\ti fixed point in $\overline X$
and, in fact, one of the following holds:
\begin{enumerate}
\item[(i)] $a$ is elliptic and $X^P\neq \varnothing$.
\item[(ii)] $a$ is parabolic and $|(\p X)^P| = 1$.
\item[(iii)] $a$ is hyperbolic and the attracting fixed point of $a$ is $P$\ti fixed.
\end{enumerate}
\end{lemma}

\begin{proof}
For~(i), we have for all $x\in X$
$$\sup_{g\in N} d(g x, x) = \sup_{k\in K_0, n\geq 0} d(k a^{-n} x, a^{-n} x).$$
Since by assumption $a$ has bounded orbits and $K_0$ is compact, the
latter quantity is bounded and hence $N$ has bounded orbits; 
the setwise decomposition $P =\langle a\rangle\cdot N$ implies that 
$P$ itself has bounded orbits
and hence admits a fixed point in $X$.

In cases~(ii) and~(iii), the set $\{a^n x : n\geq 0\}$ is unbounded for
all $x\in X$. Pick $\xi\in\p X$ and a subsequence $\{n_k\}$ of $\NN$
with $\lim_{k\to \infty} a^{n_k} x = \xi$ for all $x\in X$ (see
Proposition~\ref{prop_karlsson}). It is enough to show that $\xi$ is
$P$\ti fixed. Setting $F_j := \overline{X}^{K_j}$, we have
\begin{align*}
F_j \cap X \neq \varnothing & \kern1cm\forall\,j\in\ZZ\\
F_j \supseteq F_\ell & \kern1cm \forall\,j\leq \ell\\
a^n F_j = F_{j+n} & \kern1cm \forall\,j,n\in\ZZ\,.
\end{align*}
Picking $x\in F_j$, we have $a^{n_k} x\in F_{j+n_k} \se F_j$ for all
$k\geq 0$ and hence $\xi = \lim_{k\to\infty} a^{n_k} x \in F_j$ since
$F_j$ is closed in $\overline{X}$. Thus $\xi\in \bigcap_{j\in\ZZ} F_j
= \overline{X}^N$. It follows that $\xi$ is indeed $P$\ti fixed.
\end{proof}

We now turn back to the particular case where $X=\HH$ is the
hyperbolic space attached to a strongly nondegenerate quadratic space
$(\H,Q)$ of signature $(\alpha, 1)$. We will focus on the study of
continuous elementary representations $\pi: P \to \OO_{L_+}(Q)$ for
which $\pi(a)$ is hyperbolic with attracting fixed point $L_+$; that
is, $|\chi(a)| >1$.

\begin{prop}
\label{prop_determine_h1}%
Let $\chi: P\to \RR^*$ be a continuous homomorphism with $|\chi(a)|>1$, 
let $\ro:P\to \OO(E)$ be a continuous orthogonal representation and 
$\tau:=\chi\otimes \ro$.
\begin{enumerate}
\item[(i)] The orthogonal complement $E^{K_{-1}} \ominus E^{K_0}$ 
of $E^{K_0}$ in $E^{K_{-1}}$ is isomorphic to $\hc^1(P,\tau)$, with
isomorphism given by $v\mapsto f_{\sigma(v)}$, where for $v\in
E^{K_{-1}} \ominus E^{K_0}$, $$\sigma(v):=\sum_{n\leq -1}
\tau(a)^{n+1} v$$ and $f_{\sigma(v)}$ is the standard cocycle uniquely
determined by
$$f_{\sigma(v)}(k)=\tau(k)\sigma(v)-\sigma(v)\,.\kern1cm (\forall\,k\in K_0)$$ 

\item[(ii)] For $\eta\in\hc^1(P,\tau)$, let $f_{\sigma(v)}$ be the standard
  cocycle representing $\eta$ (with $v\in E^{K_{-1}} \ominus
  E^{K_0}$). Then the subspace $I_\eta= \overline{\langle
    f_{\sigma(v)} (p) : p\in P\rangle}$ coincides with the closed cyclic subspace
  generated by $v$ and, in fact,
\begin{equation*}
I_\eta = \overline{\langle f_{\sigma(v)}(n) : n\in N\setminus \bigcap_{j\in\ZZ} K_j\rangle}\,.
\end{equation*}
\end{enumerate}
\end{prop}

\begin{proof}
{\it (i)}~The proof consists of two steps.

{\it Claim 1:}  There is an isomorphism of topological vector spaces

\begin{align*}
\big\{v\in E : (\id - \tau(a)) v\in E^{K_0} \big\} \big/ E^{K_0}  &\lra \hc^1(P,\tau)\\
\kern1cm v &\longmapsto [f_v]\\
\end{align*}
where $f_v$ is the standard cocycle determined uniquely by
$$f_v(k) := \tau(k) v - v \kern1cm\forall\,k\in K_0.$$

Though this follows immediately from the Mayer-Vietoris sequence associated to topological HNN-extensions in continuous cohomology, we give an explicit proof:

Let $f$ be a standard cocycle.  Then 
repeated applications of the cocycle identity imply that 
for all $n\in\ZZ$ and $g\in P$ one has
\begin{equation}\label{1}
f(a^nga^{-n})=\tau(a)f(g)\,,
\end{equation}
which, applied to $n\geq0$ and $g\in K_0$ shows that 
$f$ is determined by its restriction to $K_0$.
Since $K_0$ is compact, there is $v\in E$, uniquely determined
modulo $E^{K_0}$, such that, for all $k\in K_0$,
\begin{equation}\label{2}
f(k)=\tau(k)v-v\,.
\end{equation}
Substituting (\ref{2}) into (\ref{1}) with $n=-1$ and $g\in K_0$,
we get
\begin{equation}\label{3}
(\id - \tau(a)) v\in E^{K_0}\,.
\end{equation}
Conversely, any cocycle $f$ on $K_0$ defined by (\ref{2})
with $v$ satisfying (\ref{3}) extends uniquely to a standard cocycle on $P$,
which hence completes the proof of Claim 1.

\medskip
{\it Claim 2:} The map
\begin{align*}
\sigma: E^{K_{-1}} \ominus E^{K_0}  &\lra \big\{v\in E : (\id - \tau(a)) v\in E^{K_0} \big\} \big/ E^{K_0}\\
v &\longmapsto \sum_{n\leq -1} \tau(a)^{n+1} v \kern1cm \text{mod } E^{K_0}
\end{align*}
is an isomorphism of topological vector spaces.

To start the proof of the Claim, observe that, 
since $|\chi(a)|>1$, the operator $S:= \id - \tau(a) =
\tau(a)(\tau(a)^{-1} -\id)$ has a (bounded) inverse 
given by $S^{-1} =-\sum_{n=0}^\infty \tau(a)^{-(n+1)}\,.$
Since $\tau(a)^\ell E^{K_0} = E^{K_\ell}$, this implies that 
\begin{equation}
\label{eq_j_def}%
J:=  \overline{\bigcup_{\ell\leq -1} E^{K_\ell}} \supseteq 
   \big\{v\in E : S v \in E^{K_0} \big\}\,.
\end{equation}
Defining $E_j:= E^{K_j}
\ominus E^{K_{j+1}}$, we have an orthogonal decomposition
\begin{equation}
\label{eq_j_prop}%
J = E^{K_0} \oplus \bigoplus^\wedge_{j\leq -1} E_j.
\end{equation}
Let $v = v_0 + \sum_{j\leq -1} v_j$ with $S v \in E^{K_0}$; since we
need to determine $v$ mod $E^{K_0}$, we may assume $v_0 = 0$. Then
$$S v = -\tau(a) v_{-1} - \sum_{j=1}^\infty (\tau(a) v_{-(j+1)} - v_{-j})$$
where $\tau(a) v_{-1} \in E^{K_0}$ and $\tau(a) v_{-(j+1)} - v_{-j}
\in E_{-j}$. In view of (\ref{eq_j_def}) and (\ref{eq_j_prop}), saying
that $S v\in E^{K_0}$ is equivalent to saying that
$$\tau(a) v_{-(j+1)} = v_{-j}, \kern1cm j\geq 1$$
which implies $v = \sum_{n\leq -1} \tau(a)^{n+1} v_{-1}$ and hence yields
Claim 2.

{\it (ii)}~It is clear that $I_\eta$ is contained in $\overline{\langle
  \tau(p) v : p\in P\rangle}$. Conversely, since for $k\in K_0$ we
have $f_{\sigma(v)}(k) = \tau(k) \sigma(v) -\sigma(v)$, it follows
that
$$- \int_{K_0} f_{\sigma(v)}(k) \,dk = \sigma(v) - \proj_{E^{K_0}}
(\sigma(v)) = \sigma(v)$$
and
$$ - \int_{K_{-1}} \tau(u) \int_{K_0} f_{\sigma(v)}(k) \,dk\,du = v.$$
Thus $v\in I_\eta$, which is a closed invariant subspace, and hence
contains the closure of its orbit
$\overline{\langle \tau(p) v : p\in P\rangle}$. The additional formula for $I_\eta$ follows since $f_{\sigma(v)}$ vanishes on $M=\bigcap_{j\in \ZZ}K_j$ by Lemma~\ref{lemma_std_cocycle}(ii).
\end{proof}

\section{Representations of Certain Parabolic Subgroups}
\label{sec_2_trans_par}

\subsection{Gelfand Pairs}\label{subsec_gelfand}
We start this section recalling some definitions and facts about
Gelfand pairs which will be essential in the sequel.  We refer to~\cite[Sec.~24]{Simonnet}, for example,
for a complete discussion and proofs. Whilst this theory is generally presented for unitary representations, it carries over without changes to orthogonal representations; this will be our viewpoint here.

\medskip 
Let $G$ be a locally compact group, $K<G$ a compact subgroup
and let ${\mathrm C}_{\mathrm c}(G)^{\natural_K}$ be the convolution
algebra of bi\ti $K$\ti invariant functions on $G$ with compact
support.  Then $(G,K)$ is a Gelfand pair if ${\mathrm C}_{\mathrm
  c}(G)^{\natural_K}$ is commutative.  It is easy to see that the
condition $x^{-1}\in KxK$ for all $x\in G$ is sufficient for $(G,K)$
to be a Gelfand pair.

If $(G,K)$ is a Gelfand pair, a continuous bi\ti $K$\ti invariant function 
$\varphi\in{\mathrm C}(G)^{\natural_K}$ is a spherical function if
\begin{enumerate}
\item $\varphi(e)=1$, and
\item for all $f\in{\mathrm C}_{\mathrm c}(G)^{\natural_K}$ there exists
a constant $c_f$ such that $\varphi\ast f=c_f\varphi$.
\end{enumerate}
An irreducible orthogonal representation of $G$ is $K$\ti spherical if
there exists a nonzero $K$\ti invariant vector.  If $(G,K)$ is a
Gelfand pair and $(\pi,\H)$ is any irreducible orthogonal representation
of $G$, then $\dim \H^K\leq1$, and hence $(\pi,\H)$ is $K$\ti
spherical if and only if $\dim\H^K=1$. Moreover, (equivalence classes
of) $K$\ti spherical representations of a Gelfand pair $(G,K)$ are in
bijective correspondence with positive definite spherical functions,
with the correspondence given by $\varphi(g)=\langle
\pi(g)v,v\rangle$, where $v\in\H$ is a $K$\ti invariant vector of
norm one and $\langle\,\,,\,\rangle$ is the inner product in $\H$.

\subsection{}
In the remainder of this section, we shall consider a closed subgroup
$P<\aut(\T)$ of the automorphism group of a locally finite tree $\T$
satisfying the following conditions:
\begin{enumerate} 
\item $P$ fixes a point $\xi\in \p \T$;
\item $P$ acts doubly transitively on $\p\T \setminus\{\xi\}$.
\end{enumerate}

We shall assume throughout that the vertices of $\T$ have valence at
least three; observe that under these hypotheses it follows that $P$
acts transitively on the vertices of $\T$, which is therefore a
regular tree (and $\xi$ is uniquely determined). Fix a geodesic line
$c:\ZZ\to\T$ with $c(+\infty) = \xi$; by~(2) there is a hyperbolic
element $a\in P$ with axis $c$, translation length one and attracting
fixed point $\xi$. We denote by $K_j$ the stabilizer of $c(j)$ in $P$
for $j\in \ZZ$. Then $P$ has the structure $P = \langle a\rangle \ltimes
N$ as in Section~\ref{sec_HNN}.

\bigskip

Let now $(\H, Q)$ be strongly nondegenerate of signature $(\aleph_0,
1)$ and let $L$ be an isotropic line. The main objective of this
section is to prove:

\begin{thm}
\label{thm_P_classify}%
For every continuous homomorphism $\chi: P\to \RR^*$ with
$|\chi(a)|>1$, there is up to conjugation a unique continuous
representation $\lambda: P\to\OO_L(Q)$ such that

\begin{enumerate}
\item[(i)] $P$ acts on $L$ by multiplication by $\chi$.
\item[(ii)] There is no proper closed $P$\ti invariant subspace of $\H$
  which is nondegenerate of index one.
\end{enumerate}
\end{thm}

This result is based on Section~\ref{sec_el_actions},
Section~\ref{sec_HNN} and on the following generalization of
\cite{Nebbia}:

\begin{prop}
\label{prop_rep_orth}%
Let $P_x$ be the stabilizer in $P$ of a vertex $x\in\T$. Then there is
a unique (up to equivalence) irreducible orthogonal representation of
$P$ having a $P_x$\ti fixed vector and whose restriction to $N$ is
nontrivial.
\end{prop}

We precede the proof with some intermediate results.  
For $j\in\ZZ$, we denote by $H_j$ the horosphere centered at $\xi$
passing through $c(j)$; furthermore, for $\ell\geq 0$, let
$H_j(2\ell)$ be the intersection of $H_j$ with the sphere or radius
$2\ell$ centered at $c(j)$. Notice that for $r\leq j$, the group $K_r$
acts on $H_j(2\ell)$ for every $\ell>0$, and moreover:

\begin{lemma}
\label{lemma_horos_trans}%
The group $K_r$ acts transitively on $H_j(2\ell)$ for all $r\leq j$ and $\ell\geq 0$.
\end{lemma}

\begin{proof}
  Write $\xi_\pm:=c(\pm\infty)$ and pick $a,b\in H_j(2\ell)$. Complete
  the geodesic segments $[c(j+\ell), a]$ and $[c(j+\ell), b]$ to
  infinite rays $[c(j+\ell), \alpha]$ and $[c(j+\ell), \beta]$
  respectively, where $\alpha,\beta\in\p \T$. We may assume $\ell\neq
  0$, thus $\alpha,\beta\neq\xi_-$. By double transitivity, there is
  $g\in P$ with $g(\alpha) = \beta$ and $g(\xi_-) = \xi_-$. The center
  of the tripods $(\xi_+, \alpha, \xi_-)$ and $(\xi_+, \beta, \xi_-)$
  is $c(j+\ell)$, hence $g$ fixes that point. But since $g$ fixes
  $\xi_\pm$ and hence preserves $c(\ZZ)$, it follows that
  $g\in\bigcap_{n\in\ZZ} K_n$.
\end{proof}

\begin{cor}
\label{cor_gelfand}%
$(N,K_j)$ is a Gelfand pair for all $j\in\ZZ$.\hfill\qedsymbol
\end{cor}

\begin{proof}
  As mentioned Section~\ref{subsec_gelfand}, 
  it is enough to prove that $n^{-1}\in K_jnK_j$ for all $n\in N$
  and $j\in \ZZ$.  We may assume that $n\notin K_j$. Then $n\in K_\ell
  \setminus K_{\ell-1}$ for some $\ell \geq j+1$ and $n(c(j)),
  n^{-1}(c(j)) \in H_j(2(\ell-j))$. By Lemma~\ref{lemma_horos_trans},
  this implies the existence of $k\in K_j$ with $k n(c(j)) =
  n^{-1}(c(j))$ and hence $n^{-1}\in K_jnK_j$.
\end{proof}

Let $\widehat{N}$ be the set of (equivalence classes of) irreducible
orthogonal representations of $N$; the group $\langle a\rangle$ acts on
$\widehat{N}$ by $a_*\pi(g) = \pi(a g a^{-1})$, thus preserving the
subset $\widehat{N}^1$ of representations that have a $K_j$\ti fixed
vector for some $j\in\ZZ$. In fact, if we set
$$\widehat{N}^1_j := \big\{ (\pi, \H) \in \widehat{N}^1 : \H ^{K_j} \neq 0, 
     \H^{K_{j+1}} = 0 \big\}\,,$$
     then $\widehat{N}^1 = \bigsqcup_{j\in\ZZ} \widehat{N}^1_j$ and
     $a_* \widehat{N}^1_j = \widehat{N}^1_{j-1}$, so that any
     $\widehat N_j^1$ is a fundamental domain for the action of
     $\langle a\rangle$ on $\widehat N^1$.

\begin{lemma}\label{lem_card_one}
\label{lemma_dim_one}%
$|\widehat{N}^1_j|=1$ for all $j\in\ZZ$.
\end{lemma}

\begin{proof} Since $K_\ell\leq K_j$ for $\ell\leq j$, for 
any $K_j$\ti spherical representation $(\pi,\H)$  we have
that $\H^{K_\ell}\supset\H^{K_j}$;
since $(N,K_j)$ is a Gelfand pair, these spaces are of dimension one
and hence $\H^{K_\ell}=\H^{K_j}$. Thus, in order to show that
$|\widehat{N}^1_j|=1$, it is sufficient to show that there is 
a unique positive definite $K_j$\ti spherical function $\varphi$ with 
\begin{equation*}
\int_{K_{j+1}}\varphi(kg)dk=0,\kern1cm\forall g\in N\,.
\end{equation*}
Since it suffices to show that 
$|\widehat{N}^1_0|=1$, we start  by showing that the space
$$S_0 := \Big\{\varphi\in (\cont(N)^\natural)^{K_0} : \int_{K_1}
\varphi(kg)\,dk =0 \ \forall\,g\in N \Big\}$$
is of dimension one. By identifying $N/K_0$ with the horosphere $H_0$,
we can identify $S_0$ with the space of $K_0$\ti invariant functions
on $H_0=\bigsqcup_{\ell\geq0} H_0(2\ell)$, and by applying 
Lemma~\ref{lemma_horos_trans} we deduce that any function $\varphi\in S_0$
can be written as 
\begin{equation*}
\varphi = \sum_{\ell =0}^\infty \kappa_\ell\one_{H_0(2\ell)}\,,
\end{equation*}
where $\kappa_\ell\in\RR$. The condition defining
$S_0$ means that the sum of the values of $\varphi$ over any $K_1$\ti
orbit in $H_0$ is zero; denoting by $q$ the valence of $\T$, this
implies that $\kappa_\ell =0$ for all $\ell\geq 2$ and 
$\kappa_0 +(q-2)\kappa_1 = 0$, thus proving the claim.

Thus, let $\varphi_0\in S_0$ be the unique function  such that
$\varphi(e)=1$.  To complete the proof we need to show
that for all $f\in{\mathrm C}_{\mathrm c}(G)^{\natural_K}$ there exists
a constant $c_f$ such that $\varphi\ast f=c_f\varphi$.
To this end, observe that any $f\in{\mathrm C}_{\mathrm c}(G)^{\natural_K}$
is a linear combination of characteristic functions 
$\one_{K_0}$ and $\one_{K_n \setminus K_{n-1}}:=\chi_n$ for $n>0$.
In this notation we have that
\begin{equation*}
\varphi_0=
 \one_{K_0}-\frac{1}{q-2}\one_{K_1\setminus K_0}\,.
\end{equation*}
Moreover, a direct computation shows that
$$\chi_n\ast\chi_m=
\begin{cases}
\mu(K_n\setminus K_{n-1})\chi_m&\text{if } n<m\\
\mu(K_n\setminus K_{n-1})\one_{K_n}-\mu(K_{n-1})\chi_n&\text{if }n=m\,,
\end{cases}
$$
wherein $\mu$ is the Haar measure implicit in the chosen convolution structure on ${\mathrm C}_{\mathrm c}(G)$.
Now it follows that
\begin{equation*}
\begin{aligned}
&\varphi_0\ast\chi_m\,\,=\hphantom{-}0\kern2cm\text{ if } m>1\\
&\varphi_0\ast\chi_1\,\,\,\,=-\mu(K_0)\varphi_0\\
&\varphi_0\ast\one_{K_0}=\hphantom{-}\mu(K_0)\varphi_0\,,
\end{aligned}
\end{equation*}
and hence shows that $\varphi_0$ is spherical.

Finally, observe that since $\varphi_0$ is compactly supported,
it follows from the identity 
\begin{equation*}
\varphi_0\ast\check\varphi_0=c_{\check\varphi_0}\varphi_0\,,
\end{equation*}
where as usual $\check\varphi_0(x):=\varphi_0(x^{-1})$, 
that $\varphi_0$ is positive definite.
\end{proof}

\begin{proof}[Proof of Proposition~\ref{prop_rep_orth}]
  The is no loss of generality in assuming $x=c(0)$, so that
  $P_x=K_0$.  Let $\pi_j$ be a representative of the unique
  equivalence class in $\widehat N_j^1$ and observe that
  $a_*\pi_j\cong\pi_{j-1}$.  Thus the $N$\ti representation
  $\pi:=\oplus_{j\in\ZZ}\pi_j$ extends canonically to a
  $P$-representation.  To verify that $\pi$ is irreducible, observe
  that if $\sigma$ is any sub\ti $P$\ti representation of $\pi$, then
  $\sigma|_N$ is a direct sum of sub\ti $N$\ti representations
  $\sigma_j$ of $\pi_j$, because the $\pi_j$ are irreducible and
  pairwise inequivalent.  Therefore, each $\sigma_j$ is either zero or
  irreducible.  The $a$\ti invariance of $\sigma$ shows that either
  $\sigma$ is zero or it coincides with $\pi$, hence $\pi$ is
  irreducible as a $P$\ti representation.  The uniqueness of $\pi$
  follows from Lemma~\ref{lem_card_one} and the existence of a
  $K_0$\ti fixed vector by construction.
\end{proof}

\begin{proof}[Proof of Theorem~\ref{thm_P_classify}]
  {\it Existence: } Let $\pi:P\to\OO(E)$ be the continuous
  orthogonal representation of $P$ constructed in
  Proposition~\ref{prop_rep_orth} with underlying Hilbert space $E$.
  Let $\ro:=|\chi|\otimes\pi\otimes\chi^{-1}$ and endow $E$ with the
  $P$\ti action defined by $\chi\otimes\ro$.  Since we have by
  construction that $E^{K_{-1}}\ominus E^{K_0}\neq0$, then
  Proposition~\ref{prop_determine_h1}(i) implies that
  $\hc^1(P,\chi\otimes\ro)\neq0$ and hence
  Proposition~\ref{prop_reps_isom} (with $G=P$) gives us a
  representation $\lambda:P\to\OO_L(Q)\se\OO(Q)$, so that
  $(\chi\otimes\ro)|_N=\pi|_N$ acting on $L$ by $\chi$.  By
  Proposition~\ref{prop_actions_alt}(iii), one can extract the
  corresponding ``irreducible'' part and hence the existence is
  proved.

\noindent
{\it Uniqueness:} Let us set $L_+=L$.  Since $|\chi(a)|\neq1$,
$\lambda(a)$ is hyperbolic; let $L_-$ be the isotropic line
representing the repelling fixed point and $\ro: P\to\OO(E)$ be the
orthogonal representation obtained as in Section~\ref{sec_elem}
\emph{via} the identification $E\to \pe{L_+}/L_+$ and $\eta\in\hc^1(P,
\chi\otimes\ro)$ the cohomology class defined by the above action. The
irreducibility hypothesis, Propositions~\ref{prop_dec_chi_neq_1}
and~\ref{prop_determine_h1} then imply that
\begin{equation*}
E = \overline{\langle  f_{\sigma(v)}(p) : p\in P\rangle}\,,
\end{equation*}
where $f_{\sigma(v)}$ is the standard cocycle representing $\eta$ and
moreover $E$ is the cyclic subspace generated by $v\in E^{K_{-1}}
\ominus E^{K_0}$.

Consider the orthogonal representation $\psi:=|\chi|^{-1}\otimes \chi
\otimes \ro$ on $E$ and let $\pi$ denote the orthogonal representation
of $P$ given by Proposition~\ref{prop_rep_orth}. Then
$$\psi = m\pi \oplus \pi_1 \oplus \pi_2,$$
where $\pi_1|_N = \one$ and $\pi_2$ does not have any nonzero
$K_n$\ti invariant vectors. Since $v\in E^{K_{-1}} \ominus E^{K_0}$,
the projections of $v$ to the components of $\pi_1$ and $\pi_2$ is
zero. Being a cyclic vector, this implies $\psi = m\pi$ and therefore
also $m=1$. It now follows that $\ro = |\chi| \otimes \chi^{-1}
\otimes \pi$; in the notation of Proposition~\ref{prop_determine_h1},
$\tau = |\chi|\otimes \pi$. By Proposition~\ref{prop_determine_h1}, we
have $\hc^1(P, \tau) \cong E^{K_{-1}} \ominus E^{K_0}$. Let now
$(\pi_n,\H_n)$ be the unique $K_n$\ti spherical representation of $N$
that is not $K_{n+1}$\ti spherical. Then it follows from
Lemma~\ref{lemma_dim_one} that $\tau|_N = \pi|_N = \oplus_{n\in\ZZ}
\pi_n$. Hence,
$$E^{K_{-1}} = \bigoplus_{n\geq -1}^\wedge \H_n^{K{-1}}, \kern1cm E^{K_0} = 
  \bigoplus_{n\geq 0}^\wedge \H_n^{K_0}.$$
  Observing that $ \H_n^{K{-1}} = \H_n^{K_0}$ for all $n\geq 0$, we
  deduce that $E^{K_{-1}} \ominus E^{K_0}$ has dimension one. Hence
  $\dim\hc^1(P, \tau) = 1$, which implies now by
  Proposition~\ref{prop_reps_isom} that, up to conjugation, $\lambda:
  P\to \OO_L(Q)$ is completely determined by $\chi$.
\end{proof}

\section{Representations of $G$ into $\OO(Q)$}
\label{sec_rep}

\subsection{}
In this section $\T$ denotes a regular or biregular tree of
finite bivalency $(r,s)$ with $r,s\geq3$.  A subgroup
$G<\aut(\T)$ satisfies the property $T^+_2$ if for every
$\xi_1\neq\xi_2$ in $\p\T$ and
$\eta_1,\eta_2\in\p\T\setminus\{\xi_1,\xi_2\}$ such that
the distance 
between the projections of $\eta_1$ and $\eta_2$ 
on the geodesic $[\xi_1,\xi_2]$ is even, 
there exists $h\in G$ fixing $\xi_1,\xi_2$ and
$h(\eta_1)=\eta_2$.

This property is implied by triple transitivity of the $G$\ti action on
$\p\T$ and implies double transitivity.  The main
result of this section is:

\begin{thm} 
\label{thm_indef_comp}%
Let $G<\aut(\T)$ be a closed subgroup satisfying property $T^+_2$,
$\xi\in\p\T$ and $P$ the stabilizer of $\xi$ in $G$.  Let $(\H,Q)$ be
a strongly nondegenerate quadratic space of index 1 and
$\pi:G\to\OO(Q)$ a continuous, nonelementary representation.  Then
$\pi|_P$ has an irreducible indefinite component $\H_{1,P}$ and the
canonical orthogonal decomposition $\H=\H_{1,P}\oplus\H_{0,P}$ is
$G$\ti invariant.

Let $\pi_1,\pi_2:G\to\OO(Q)$ be nonelementary continuous
representations such that $\pi_1|_P=\pi_2|_P$.  Then the restriction
of $\pi_1$ and $\pi_2$ to the indefinite irreducible components
of $P$ coincide.
\end{thm}

\subsection{}
\label{sec_bnd_gen}%
Let $G$ be a locally compact group boundedly generated by $\{s\}\cup P$, where
$s\in G$ and $P<G$ is a closed subgroup with the structure considered in Section~\ref{sec_HNN}; assume further that:
\begin{enumerate}
\item $\langle s\rangle$ is relatively compact, and
\item $\{a^nsa^ns^{-1}:\,n\geq1\}$ is relatively compact. 
\end{enumerate}

\noindent
(We recall that a group $G$ is said \emph{boundedly generated} by a subset $X\se G$ if there is $n\in\NN$ with $X^n = G$.)

\begin{prop}\label{prop_7_2} 
  Let $X$ be a complete \cat{-1} space and $G\times X\to X$ a
  continuous, nonelementary isometric action.  Then $a$ acts as a
  hyperbolic element, $s$ exchanges both fixed points of $a$ and $(\p
  X)^P$ is the attracting fixed point of $a$.
\end{prop}

\begin{proof} 
  We use Lemma~\ref{lem_class_isom} and we distinguish the three cases for $a$.

\smallskip

{\it $a$ is elliptic}: Then $X^P\neq\emptyset$ and since $\langle
s\rangle$ is relatively compact and $G$ is boundedly generated by
$\langle s\rangle$ and $P$, the $G$\ti orbits in $X$ are bounded.
Hence $X^G\neq\emptyset$, contradicting non-elementarity.

\smallskip

{\it $a$ is parabolic}: Then $\p X^P=\{\xi\}$.  Let
$\{n_i\}_{i\geq1}$ be a sequence such that $a^{n_i}x\to\xi$; then
$a^{-n_i}x\to\xi$, and since $\{a^nsa^ns^{-1}:n\geq1\}$ is bounded, we
have that $sa^{n_i}s^{-1}x\to\xi$ which, in view of the fact that
$a^{n_i}s^{-1}x\to\xi$, implies that $s(\xi)=\xi$ and hence
$G\xi=\xi$, contradiction.

\smallskip

Thus, {\it $a$ is hyperbolic}. Let now $\xi_-,\xi_+$ be respectively the
repelling and the attracting fixed points of $a$ on $\p X$.  Since
$a^nx\to\xi_+$, $a^{-n}x\to\xi_-$, and $\{a^nsa^ns^{-1}:n\geq1\}$ is
bounded, we deduce that $sa^ns^{-1}x\to\xi_-$ and hence
$s(\xi_+)=\xi_-$ and $s(\xi_-)=\xi_+$.

Finally, we know $(\p X)^P\ni \xi_+$ from Lemma~\ref{lem_class_isom}; since $(\p X)^{\langle a\rangle} = \{\xi_\pm\}$, it remains only to observe that $\xi_-$ is not $P$\ti fixed, since otherwise the set $\{\xi_\pm\}$ would be preserved by $G$.
\end{proof}

\subsection{}
\label{subsec_snshns}%
Let now $G<\aut(\T)$ be any closed subgroup which acts doubly transitively
on $\p\T$. Then (see~\cite{Burger_Mozes_IHES_2} Sec.~4.1 and~0.4):

\begin{itemize}
\item[--] For every $\xi\in\p\T$, the Busemann character
$\chi_\xi:G_\xi\to\ZZ$ has image $\ZZ$ or $2\ZZ$ depending on whether
$G$ is vertex transitive or not;
\item[--] For every geodesic $c:\ZZ\to\T$ there is $s\in
G$ and $n_0\in\ZZ$ with $s(c(\pm\infty))=c(\mp\infty)$ and $s(c(n_0))=c(n_0)$;
\item[--] For any $\xi\neq\eta$ in $\p\T$ and $s\in G$
exchanging $\xi$ and $\eta$, we have $G=G_\xi\cup G_\xi sG_\xi$.
\end{itemize}

\begin{lemma}\label{lem_snshns}
  Let $G<\aut(\T)$ be a closed subgroup satisfying $T^+_2$,
  $c:\ZZ\to\T$ a geodesic, $\xi_\pm=c(\pm\infty)$, and $S=\{s\in
  G:s(\xi_\pm)=(\xi_\mp),\,s(c(0))=c(0)\}$.  Let also $K_j=G_\xi\cap
  G(c(j))$.  Then, given $j\in\ZZ$, $n\in K_j\setminus K_{j-1}$ and
  $s\in S$, there exists $h\in G$ such that:
\begin{enumerate}
\item[(i)] For all $q\in \ZZ$, $hc(q)=c(q-2j)$;
\item[(ii)] For all $s',s''\in S$, $s'nshns''\in K_{-j}\setminus K_{-j-1}$. 
\end{enumerate}
\end{lemma}

\begin{proof}
  Since $n\in K_j\setminus K_{j-1}$, then
  $\proj_{[\xi_+,\xi_-]}(n\xi_-)=c(j)$ and
  $\proj_{[\xi_+,\xi_-]}(s^{-1}n^{-1}\xi_-)=c(-j)$.  Thus property
  $T_2^+$ implies that there exists $h\in G_{\xi_+}\cap G_{\xi_-}$
  such that $hn\xi_-=s^{-1}n^{-1}\xi_-$ so that (i) follows.  Let now
  $s',s''\in S$ and set $g=s'nshns''$.  Then we have:

\smallskip

($\circ$)~$g\xi_+=s'nshns''\xi_+=s'nshn\xi_-=s'nss^{-1}n^{-1}\xi_-=s'\xi_-=\xi_+$.

\smallskip

($\circ$)~$gc(-j)=s'nshns''c(-j)=s'nshnc(j)=s'nshc(j)=s'nsc(-j)=s'nc(j)=s'c(j)=c(-j)$.
Thus $g\in K_{-j}$.

\smallskip

($\circ$)~$gc(-j-1)=s'nshns''c(-j-1)=s'nshnc(j+1)=s'nshc(j+1)=s'nsc(-j+1)=s'nc(j-1)$.
But since $n\in K_j$, we have that $nc(j-1)\neq c(j+1)$ and hence
$s'nc(j-1)\neq c(-j-1)$. Thus $g\notin K_{-j-1}$.
\end{proof}

\subsection{}

We are finally ready to give the 

\begin{proof}[Proof of  Theorem~\ref{thm_indef_comp}] 
Let $\pi:G\to\OO(Q)\to\isom(\HH)$ be a continuous nonelementary action. Being the stabilizer of $\xi\in\p\T$, $P$ has the structure of Section~\ref{sec_HNN}; we shall use the corresponding notation for $a$. In view of Section~\ref{subsec_snshns}, we can choose $s$ and a parametrisation $c:\ZZ\to\T$ of the axis of $a$ such that $s c(0) = c(0)$ and $s c(\pm\infty) = c(\mp\infty)$. Observe further that these notations also put us in the setting of Section~\ref{sec_bnd_gen}. By Proposition~\ref{prop_7_2}, $\pi(a)$ is hyperbolic and hence cannot fix a point in $\HH$ or preserve any horosphere. It follows from Proposition~\ref{prop_actions_alt} applied to $P$ that $\pi|_P$ has an irreducible indefinite component $\H_{1,P}$.
  
  We need to show that $\H_{1,P}$ is $G$\ti invariant and that, on
  $\H_{1,P}$, the representation $\pi$ is determined by $\pi|_P$.  Let
  $L_\pm$ be the attr./repell. fixed points of $\pi(a)$ and let
  $L=L_+\oplus L_-$, $F=\pe L$.  We are in the setting of
  Proposition~\ref{prop_actions_alt} for $P$ (instead of $G$) 
  and we adopt its notation.  By Proposition~\ref{prop_dec_chi_neq_1}, 
  $\H_{1,P}=L\oplus j^{-1}(I_\eta)$. By Proposition~\ref{prop_7_2}, 
  $\pi(s)$ exchanges $L_\pm$.  It is enough to show that $\pi(s)$ preserves
  $j^{-1}(I_\eta)$ and that its restriction to $\H_{1,P}$ is determined
  by $\pi|_P$.
  
  We adopt the notation of Appendix~\ref{sec_app} with $\ell_\pm$.
  Then
\begin{equation*}
\pi(s)=
\begin{pmatrix}
0&\mu&0\\
\mu^{-1}&0&0\\
0&0&\pi_0(s)
\end{pmatrix}\,,
\end{equation*}
where $\pi_0(s)$ is orthogonal and $\mu\in\RR^*$.

Fix $j\in\ZZ$, $n\in K_j\setminus K_{j-1}$ and $h$ as in
Lemma~\ref{lem_snshns} such that $g:=snshns\in K_{-j}\setminus
K_{-j-1}$.  We write
\begin{equation*}
\pi(n)=
\begin{pmatrix}
\chi(n)&\alpha(n)&N_2^+\\
0&\chi(n)^{-1}&0\\
0&N_3^-&\pi_0(n)
\end{pmatrix}\,,
\end{equation*}

and likewise
\begin{equation*}
\pi(g)=
\begin{pmatrix}
\chi(g)&\alpha(g)&M_2^+\\
0&\chi(g)^{-1}&0\\
0&M_3^-&\pi_0(g)
\end{pmatrix}\,.
\end{equation*}

As to $\pi(h)$, since it fixes both $L_\pm$, it is of the form

\begin{equation*}
\pi(h)=
\begin{pmatrix}
\chi(h)&0&0\\
0&\chi(h)^{-1}&0\\
0&0&\pi_0(h)
\end{pmatrix}\,.
\end{equation*}

Computing $g=snshns$ we find
\begin{equation}\label{eq:7.1}
\chi(g)=\mu^{-1}\chi(n^{-1}h)\alpha(n)\,,
\end{equation}
and
\begin{equation}\label{eq:7.2}
M_3^-=\chi(hn)\pi_0(s)N_3^-\,.
\end{equation}

Equation (\ref{eq:7.1}) shows that $\mu$ is determined in terms of $\pi|_P$.  We
are left to determine $\pi_0(s)$.  Write $f$ for the standard cocycle
associated to $\eta$; then (\ref{eq:7.2}) gives
\begin{equation*}
\pi_0(s)f(n)=\chi(gh)^{-1}f(g)\,.
\end{equation*}
We obtain such a formula for every $n\in N\setminus\cap_{j\in\ZZ}K_j$.
Thus we are done since by Proposition~\ref{prop_determine_h1}(ii),
$I_\eta$ is spanned by these $f(n)$.
\end{proof}

\begin{proof}[Proof of Theorem~\ref{thm_class_intro}]
By Proposition~\ref{prop_isom_H}, we are reduced to study homomorphisms from
$G$ into $\OO_+(Q)$ and can thus apply the results obtained so far.
The existence statement in Theorem~\ref{thm_class_intro} follows 
from Theorem~\ref{thm_constr_intro} and Proposition~\ref{prop_actions_alt}(iii).

For the uniqueness part, let $\xi\in\partial\T_r$,
$P$ the stabilizer of $\xi$ in $G$ and $a\in P$ a hyperbolic element
with attracting fixed point $\xi$ and translation length $1$.
Let $\pi:G\to\OO_+(Q)$ be a nonelementary continuous representation.
Let $c:\ZZ\to\T_r$ be a parametrization of the axis of $a$,
and let $K_n$ be, as usual, the stabilizer in $P$ of $c(n)$.
Using that $G$ is doubly transitive on $\partial\T_r$,
we see that any other hyperbolic element $b$ with translation length $1$ 
is conjugate to an element of the form $a\cdot k$, 
where $k\in\cap_{n\in\ZZ}K_n$.  Since $\pi(a)$ is hyperbolic
(Proposition~\ref{prop_7_2}), and using that $a$ normalises $\cap_{n\in\ZZ}K_n$,
one sees that $\pi(k)$ fixes pointwise the axis of $\pi(a)$
and hence $\pi(a)$ and $\pi(b)$ have the same translation length,
say $\ell_\pi$.  Let $L_+$ be the attractive fixed point of $\pi(a)$
and $\chi$ the character by which $P$ acts on $L_+$.  
Then, since $\pi$ takes values in $\OO_+(Q)$, 
we have that $\chi$ takes values in $\RR_+$ which implies first that 
$\chi(a)=e^{\ell_\pi}$ and then that $\chi$ is trivial on $N$,
since $N$ is an increasing union of compact groups.  
This shows that $\chi$ is completely determined by $\ell_\pi$.

Assume now that $\pi$ is irreducible.  Then Theorem~\ref{thm_indef_comp} implies
that $\pi|_P$ is irreducible, Theorem~\ref{thm_P_classify} that it is completely
determined by $\chi$ (and hence by $\ell_\pi$), and Theorem~\ref{thm_indef_comp}
again that $\pi$ is completely determined by $\ell_\pi$.

\end{proof}
\section{Explicit Constructions}
\label{sec_constr}%

\subsection{}
\label{sec_constr_gen}%

Let $\T$ be any simplicial tree. Let $\alpha+1$ be the cardinal of the
vertex set $V$ of $\T$ and let $(\H, Q)$ be a strongly nondegenerate
quadratic space of signature $(\alpha, 1)$; let $\HH$ be the
corresponding hyperbolic space. Denote by $G$ the (abstract) group
$G=\aut(\T)$. We denote by $d$ both the metric on $\HH$ and the metric
on (the geometric realization of) $\T$ that gives unit length edges.

\begin{thm}
\label{thm_constr}%
For every $\lambda>1$ there is an embedding $\Psi:\T\to\HH$ and a
representation $\pi:G\to\OO(Q)\to\isom(\HH)$ such that:
\begin{enumerate}
\item[(i)] The map $\Psi$ is $G$\ti equivariant for $\pi$.
\item[(ii)] $\lambda^{d(x,y)} = \cosh d(\Psi x, \Psi y)$ for any two
  vertices $x,y$ of $\T$
\item[(iii)] $\Psi$ extends to an equivariant
  boundary map $\p\Psi: \p T\to\p\HH$ which is a homeomorphism onto
  its image.
\item[(iv)] $\Psi(V)$ is cobounded in the convex hull $\C\se \HH$ of the
  image of $\p\Psi$.
\end{enumerate}
\end{thm}

\begin{rem}
\label{rem_Busemann}%
The formula in (ii) shows in particular that $d(\Psi x, \Psi y)$ is
asymptotically proportional to $d(x,y)$. If we denote by $b^\T$ the
Busemann cocycle for $\T$ and $b$ is the one for $\HH$ as in
Section~\ref{sec_hyperbolic}, we have
$$b_{\Psi\xi}(\Psi x, \Psi y) = b^\T_\xi(x,y) \ln\lambda
\kern1cm\forall\,\xi\in\p\T, \forall\, x,y\in V.$$
\end{rem}

We can give right away the construction of $\Psi$; the remainder of
the section will be devoted to proving the properties stated in
Theorem~\ref{thm_constr}.

\medskip

Fix a vertex $w\in V$. By Proposition~\ref{prop_signature}, we may
identify $\H$ with $\ell^2(V)$ in such a way that the bilinear form
$B$ associated to $Q$ reads
$$B(f,g) = \sum_{v\in V, v\neq w} f(v)g(v) - f(w) g(w).$$
We define a map $V\to\H$, $v\mapsto f_v$ as follows. Denote for $u\in
V$ by $\delta_u$ the unit function supported on $u$; then
$$f_v:=\lambda^{d(w,v)} \delta_w + \sqrt{\lambda^2
  -1}\sum_{k=1}^{d(w,v)} \lambda^{d(w,v)-k}\,\delta_{u_k},$$
where $w, u_1, u_2, \ldots, u_{d(w,v)}=v$ is the geodesic from $w$ to
$v$ (it is understood that the right hand side summation is zero when
$v=w$). A computation gives $Q(f_v)=-1$ so that $f_v$ is in the
negative cone $C_-$; now $\Psi$ is the resulting map $V\to C_- \to
\HH$ extended to $\T$ by sending each edge to a geodesic segment. Each
element $\xi\in\p \T$ can be realized by a unique geodesic ray of
vertices $\{v_k\}_{k=0}^\infty$ with $v_0 = w$; we define $(\p\Psi)(\xi)$
by considering the element $f_\xi$ of the isotropic cone $C_0$ given
by
$$f_\xi := \delta_w + \sqrt{\lambda^2 -1} \sum_{k=1}^\infty\lambda^{-k} \delta_{v_k}.$$
Observe that one obtains a multiple of $f_{x_k}$, hence the same point
$\Psi x_k$, by truncating the above sum at $k$. It follows that the
resulting map $\overline{\T}\to\overline{\HH}$ is continuous; the
claim~(iii) now follows from Remark~\ref{rem_Busemann} and claim~(i).
The formula of claim~(ii) can be verified by inspection; however, we
shall see that it can be reduced to the obvious case $u=w$.

\bigskip

The strategy for the proof of Theorem~\ref{thm_constr} is first to
construct $\pi$ in the case where $\T$ is \emph{regular}, that is, $G$
acts transitively on $V$. The general case will follow by the
naturality of our construction with respect to the pointed tree
$(\T,w)$.

\subsection{Regular Case}

We assume that $G=\aut(\T)$ acts transitively on $V$. Fix a neighbour
$z$ of $w$ and denote by $K<G$ the stabilizer of $w$, by $L$ the
subgroup of $G$ preserving the set $\{w,z\}$ and by $E = K\cap L$ the
pointwise stabilizer of $\{w,z\}$. The $K$\ti action on $V$ preserves $B$ and hence
induces a representation $\pi: K\to \OO(Q)$ by $\pi(g)\delta_u =
\delta_{gu}$ ($g\in K$). We extend now $\pi$ by defining $\pi: L =
E\sqcup (L\setminus E)\to \OO(Q)$ as follows: the map is already
defined on $E=K\cap L$; for every $g\in L\setminus E$ and every $u\in
V$, set
\begin{equation}\label{eq_rep_2}
\pi(g)\delta_u :=
\begin{cases}
\lambda \delta_w + \sqrt{\lambda^2 -1} \delta_{z} & \text{if $u=w$.}\\
-\sqrt{\lambda^2 -1} \delta_w - \lambda \delta_{z} & \text{if $u=z$.}\\
\delta_{gu} & \text{otherwise.}
\end{cases}
\end{equation}
It is a matter of computation to verify that $\pi(g)$ is in $\OO(Q)$.

\begin{prop}
The map $\pi: L\to \OO(Q)$ extends
uniquely to a homomorphism $\pi: G\to\OO(Q)$.
\end{prop}

\begin{proof}
  We start by showing that the map $\pi:L\to\OO(Q)$ is a homomorphism;
  that is, we need to verify that $\pi(g)\pi(g') = \pi(g g')$ holds on $L$,
  which we do by discussing the cases according to where
  $g,g'$ are in the coset decomposition $L = E\sqcup (L\setminus E)$. There
  is nothing to do if $g,g'$ are both in $E$ since $\pi$ is a
  homomorphism on $K$. We shall write out the verification in the case
  $g,g'\in L\setminus E$; the two remaining cases are simpler and
  similar. Let thus $g,g'\in L\setminus E$. Then $\pi(g g')\delta_u =
  \delta_{g g' u}$ for all $u\in V$ since $g g'\in E$ (as $E$ has
  index two in $L$). On the other hand, we have:

{\it Case $u=w$}:
\begin{multline*}
  \pi(g)\pi(g')\delta_w = \pi(g)(\lambda \delta_w + \sqrt{\lambda^2 -1} \delta_{z}) =\\
  = \lambda (\lambda \delta_w + \sqrt{\lambda^2 -1} \delta_{z})
  +\sqrt{\lambda^2 -1}(-\sqrt{\lambda^2 -1} \delta_w - \lambda
  \delta_{z}) = \delta_w,
\end{multline*}
which is indeed $\delta_{g g' u}$ since $E$ fixes $w$.

{\it Case $u=z$}:
\begin{multline*}
  \pi(g)\pi(g')\delta_z = \pi(g)(-\sqrt{\lambda^2 -1} \delta_w - \lambda \delta_{z}) =\\
  = -\sqrt{\lambda^2 -1}(\lambda \delta_w + \sqrt{\lambda^2 -1}
  \delta_{z}) -\lambda(-\sqrt{\lambda^2 -1} \delta_w - \lambda
  \delta_{z}) = \delta_z,
\end{multline*}
which is indeed $\delta_{g g' u}$ since $E$ fixes $z$.

{\it Case $u\neq w,z$}: then we have also $g'u\neq w,z$ and hence
$$\pi(g)\pi(g')\delta_u = \pi(g) \delta_{g' u} = \delta_{g g' u} = \pi(g g')\delta_u.$$

  To show that the map $\pi$ defined on $L$ extends uniquely to a homomorphism 
  on $G$, observe the $G$\ti action on (the first barycentric subdivision of) $\T$
  determines a splitting of $G$ into an amalgamation $G=K \amal E L$.
  Therefore the statement follows from the universal property of
  amalgamations.
\end{proof}

\begin{rem}
\label{rem_dep_z}%
The definition of $\pi$ on $G$ is independent of the choice of $z$:
Indeed, observe first that $K$ acts transitively on the set of
neighbours of $w$. If $k$ is any element of $K$, we obtain another
neighbour $k z$ of $w$ and another amalgamation $G=K \amal{k E k^{-1}}
k Lk^{-1}$. With $\pi$ defined as before using $L$, one checks
immediately that for $g\in k L k^{-1} \setminus k E k^{-1}$ the
formula~(\ref{eq_rep_2}) for $\pi(g)$ remains valid upon replacing $z$
with $k z$.
\end{rem}

We turn to point~(i) in Theorem~\ref{thm_constr}. Pick a vertex $v\in
V$ and let $n=d(v,w)$. There is a hyperbolic element $a$ of
translation length one admitting an axis $\{u_k\}_{k\in \ZZ}$ such
that $u_0 =w$, $u_n = v$ and $a u_k = u_{k+1}$ for all $k\in \ZZ$.
Notice that $K$ and $a$ generate $G$; since moreover $\Psi$ is $K$\ti
invariant by its construction, we need only verify that $\Psi(a v) =
\pi(a)\Psi(v)$. By Remark~\ref{rem_dep_z} there is no loss of generality in
assuming that $z = u_1$. Let $s$ be an element of $L\setminus E$
preserving $\{u_k\}$; that is, $s u_k = u_{1-k}$ for all $k\in \ZZ$.
Then $a= s t$ for $t\in K$ such that $t u_k = u_{-k}$ for all $k\in
\ZZ$, and thus an immediate computation using~(\ref{eq_rep_2}) for
$\pi(s)$ shows that we have for all $k\in \ZZ$
$$\pi(a) \delta_{u_k}= \pi(s)\pi(t)\delta_{u_k}=
\begin{cases}
\lambda \delta_w + \sqrt{\lambda^2 -1} \delta_{u_1} & \text{if $k=0$.}\\
-\sqrt{\lambda^2 -1} \delta_w - \lambda \delta_{u_1} & \text{if $k=-1$.}\\
\delta_{u_{k+1}} & \text{otherwise.}
\end{cases}
$$
Now we can compute
\begin{align*}
  \pi(a) f_v & = \pi(a) \Big( \lambda^n \delta_w + \sqrt{\lambda^2 -1} \sum_{k=1}^n \lambda^{n-k} \delta_{u_k} \Big)\\
  & =\lambda^n (\lambda \delta_w + \sqrt{\lambda^2 -1} \delta_{u_1}) + \sqrt{\lambda^2 - 1} \sum_{k=1}^n \lambda^{n-k} \delta_{u_{k+1}}\\
  & = \lambda^{n+1} \delta_w + \sqrt{\lambda^2 -1} \sum_{k=1}^{n+1}
  \lambda^{n+1-k} \delta_{u_k} = f_{u_{k+1}} = f_{a v},
\end{align*}
and claim~(i) is proved. Now the transitivity of $G$ reduces~(ii) to
the case where one of the vertices is $w$, which is an immediate
computation. As~(iii) was addressed before, we are left with
proving~(iv).

\begin{prop}
  Every point $x\in \C$ is at distance at most $\cosh^{-1}
  \sqrt{1+\lambda}$ of some element of $\Psi(V)$.
\end{prop}

\begin{proof}
  Every element of $\HH$ is represented by a unique function in
  $\ell^2(V)$ with value one on $w$. In fact, if $D'$ denotes the unit
  ball in $\ell^2(V\setminus\{w\})$, the set of such functions is
  $D:=\delta_w + D'$. This gives the \emph{Klein model} in finite
  dimensions, and thus it follows that geodesics in $\HH$ correspond
  to affine lines in $D$. It is therefore enough to show the claim for
  any finite convex combination $x= \sum_{\xi\in S} s_\xi f_\xi$ where
  $S\se \p\T$ is a finite set and $s: S\to (0,1)$ is any function
  $\xi\mapsto s_\xi$ of sum one (observe that $S$ must contain at
  least two points since $0< s_\xi < 1$).
  
  We write $\beta = b(\cdot, \delta_w)$ for Busemann function on $\HH$
  normalised at $\delta_w$ and $\beta^\T = b^\T(\cdot, w)$ for the
  analogous Busemann function on $\T$. Consider the function $\psi_s:
  V \to \RR^*_+$ defined by
$$\psi_s(v) := \sum_{\xi\in S} s_\xi \lambda^{\beta^\T_\xi(v)}.$$
This function admits a minimum $v_0\in V$; indeed, outside a finite
subset of $V$ determined by the configuration of $S$, it increases
monotonically with the distance to this subset. We shall see that
$\Psi(v_0)$ is at distance at most $\cosh^{-1} \sqrt{1+\lambda}$ of
$x$.

Pick $g\in G$ such that $g v_0 = w$; then $w$ is a minimum of the
function $v\mapsto \psi_s(g^{-1}v)$, which in view of
$\beta^\T_\xi(g^{-1} v) = \beta^\T_{g \xi}(v) - \beta^\T_{g\xi}(g w)$
reads
$$\psi_s(g^{-1}v) = \sum_{\eta\in g S} s_{g^{-1}\eta}
\lambda^{\beta^\T_\eta(v)-\beta^\T_\eta(g w)}.$$
Thus, setting $(g_\star s)_\xi := \lambda^{-\beta^\T_\xi(g w)}
s_{g^{-1}\xi}$, it follows that $w$ is also a minimum of the function
$\psi_{s'}$ for $s' := g_\star s / \sum_\xi (g_\star s)_\xi$. Setting
$$\sigma_v:= \sum\big\{ s'_\xi : v \text{ is in the ray }
[w,\xi]\big\} \kern1cm \forall\, v\in V$$
 the minimality implies for every neighbour $v$ of $w$
%
 \begin{align*}
1 = \psi_{s'}(w) \leq \psi_{s'}(v) = &\lambda^{-1} \sum \{ s'_\xi :
 v\in [w,\xi]\} + \lambda \sum \{ s'_\xi : v\notin [w,\xi]\} =\\
 &\lambda^{-1} \sigma_v + \lambda (1-\sigma_v)\,,
\end{align*}
%
and hence $\sigma_v \leq \lambda / (1+\lambda)$. This in turn implies
$$\sigma_v \leq \frac{\lambda}{1+\lambda} \kern1cm \forall\, v\in V, v\neq w.$$
The formula~(\ref{eq_buse}) in Section~\ref{sec_hyperbolic} shows that for every
$g\in G$ and $h\in C_0$ we have
$$\beta_{g h} (g \delta_w) = \ln\frac{B(g h, g\delta_w)
  \sqrt{-Q(\delta_w)}}{B(g h, \delta_w)\sqrt{-Q(g\delta_w)}} =
\ln\frac{B(h,\delta_w)}{B(g h, \delta_w)} = \ln\frac{h(w)} {(g
  h)(w)}.$$
Thus we have for all $\xi\in\p\T$
$$(g f_\xi) ( w) f_\xi(w) e^{-\beta_{g f_\xi} (g \delta_w)} = f_\xi(w)
\lambda^{-\beta^\T_{g \xi}(g w)}$$
and since $g f_\xi$ is proportional to $f_{g \xi}$ we deduce $g f_\xi
= \lambda^{-\beta^\T_{g \xi}(g w)} f_{g \xi}$.  We conclude that $g x$
is represented in $\HH$ by the element $y := \sum_{\xi\in g S} s'_\xi
f_\xi$ of $D$. We proceed now to compute
$$\cosh d(\delta_w, y) = -\frac{B(\delta_w, y)}{\sqrt{Q(\delta_w)
    Q(y)}} = \frac{1}{\sqrt{-Q(y)}}.$$
If $S_k$ is the sphere of radius $k$ around $w$, we deduce from the
definition of $f_\xi$
$$y = \delta_w + \sqrt{\lambda^2 -1} \sum_{k=1}^\infty \sum_{v\in S_k}
\lambda^{-k} \sigma_v\delta_v$$
and therefore
$$Q(y) = -1 + (\lambda^2 -1) \sum_{k=1}^\infty \lambda^{-2
  k}\sum_{v\in S_k} \sigma_v^2 \leq -1 + (\lambda^2 -1)
\frac{\lambda}{1+\lambda} \sum_{k=1}^\infty \lambda^{-2 k}\sum_{v\in
  S_k} \sigma_v.$$
Using $\sum_{v\in S_k} \sigma_v = 1$ one gets finally $Q(y) \leq
-1/(1+\lambda)$. It follows that
$$d(\Psi v_0, x) = d(\Psi w, g x) = d(\delta_w, y) \leq \cosh^{-1}
\sqrt{1+\lambda}.$$
\end{proof}

This proposition completes the proof of Theorem~\ref{thm_constr} in
the regular case.\hfill\qedsymbol

\subsection{General Case}

Suppose now that $\T$ is a general tree with vertex set $V$ of
cardinal $\alpha+1$ and $G=\aut(\T)$. Complete $\T$ to a regular tree
$\T'$ with vertex set $V'\supseteq V$ of cardinal $\alpha'+1$. We keep
the notation of Section~\ref{sec_constr_gen} and define likewise $\H'
= \ell^2(V')$, $\HH'$, $B'$, $G'=\aut(\T')$, etc.; take $w':=w\in V\se
V'$ and observe that $\HH$ is a hyperbolic subspace of $\HH'$. Let
$\Psi', \pi'$ be the maps associated to $\T'$ by the proof for the
regular case.

Denote by $L_0< L< G'$ the pointwise stabilizer, respectively stabilizer, of $\T$;
since any automorphism of $\T$ can be extended to some automorphism(s)
of $\T'$, we have a natural identification $G = L/L_0$. It follows at
once from the definition of $\Psi, \Psi'$ and $\pi, \pi'$ that the
restriction of $\Psi'$ gives $\Psi$, while $\pi'$ descends to $\pi$.
In fact, as all definitions for $\T'$ vanish on $V'\setminus V$ when
restricted to $\T$, the only part of Theorem~\ref{thm_constr} for $\T$
that does not follows immediately from the case of $\T'$ is point~(iv).
But the proof given above, when applied to $\T'$ and to the
corresponding convex hull $\C'$, shows in fact that whenever $x\in
\C\se \C'$ is a finite affine convex combination of elements in
$\Psi(\p\T)$, then the vertex $v_0\in V'$ is actually in $V$. This
follows indeed from the definition of the function $\psi_s$ and thus
concludes the proof of Theorem~\ref{thm_constr}.\hfill\qedsymbol

\appendix

\section{Matrix Representations}\label{sec_app}

Let $(\H,Q)$ a strongly nondegenerate quadratic space of signature
$(\alpha,1)$.  If $L_-,L_+$ are two distinct isotropic lines, define
$L:=L_+\oplus L_-$ and $F=\pe L$, so that
$$\H=L\oplus F.$$
On $L$ and $F$ we consider the restrictions
$B|_{L\times L}$ and $B|_{F\times F}$, and if $A$ is a continuous
linear operator between any of these spaces, $A^*$ denotes the adjoint
with respect to these strongly nondegenerate bilinear forms.  For any
continuous linear operator $T:\H\to\H$, define
$$
\begin{aligned}
A_1=p_LT|_L,\qquad&A_2=P_LT|_F,\\
A_3=P_FT|_L,\qquad&A_4=p_FT|_F\,.
\end{aligned}$$
Then $T\in \OO(Q)$ if and only if the following conditions are satisfied:
\begin{enumerate}
\item[$\ $] $A^*_1A_1+A^*_3A_3=\id_L$;
\item[$\ $] $A^*_2A_2+A^*_4A_4=\id_F$\,.
\item[$\ $] $A^*_1A_2+A^*_3A_4=0$;
\item[$\ $] $A^*_2A_1+A^*_4A_3=0$;
\end{enumerate}
Observe that by taking adjoints, the last two conditions are equivalent.

We shall look more closely at $\OO_{L_+}(Q)$, the stabilizer in
$\OO(Q)$ of $L_+$.  For this, let $L_\pm=\RR\ell_\pm$, with
$B(\ell_+,\ell_-)=1$.  We represent $A_1$ by a two-by-two real matrix;
$A_2:F\to L$ will be represented by two linear forms $A_2^+$ and
$A_2^-$ given by $A_2(e)=A_2^+(e)\ell_++A_2^-(e)\ell_-$; $A_3:L\to F$
will be represented by two vectors $A_3^+=A_3(\ell_+)$,
$A_3^-=A_3(\ell_-)$, and hence
$$T=\bpm A_1&\begin{matrix} A_2^+\\A_2^-\end{matrix}\\
\begin{matrix}A_3^+&A_3^-\end{matrix}&A_4\epm\,.$$
Then $T\in\OO_{L_+}(Q)$ if and only if it has the form
\begin{equation*}
\bpm 
\lambda &       \alpha      & A_2^+ \\
   0    & \lambda^{-1} &   0   \\
   0    &    A_3^-     &   A_4
\epm
\end{equation*}
with $\lambda\in\RR^\times$, $A_3^-\in F$, $A_4\in\OO(F)$ and $\alpha,
A_2^+$ are determined by
\begin{equation*}
\alpha=-\frac\lambda2Q(A_3^-)\,,
\end{equation*}
and
\begin{equation*}
A_2^+(v)=-\lambda B(A_4(v),A_3^-), \kern1cm\text{ for all }v\in F\,.
\end{equation*}
The inverse of 
$$
S=
\bpm
\mu  &      \beta   &  B_2^+\\
 0   & \mu^{-1} &  0 \\
 0   &   B_3^-  &  B_4
\epm\in\OO_{L_+}(Q)
$$
is given by 
$$
S^{-1}=
\bpm
\mu^{-1} &         \beta           & -\mu^{-1}B_2^+B_4^{-1}\\
    0    &        \mu           &             0         \\
    0    & -\mu B_4^{-1}(B_3^-) &           B_4^{-1}
\epm
$$
and the conjugate $STS^{-1}$ of $T$ by $S$
has the following entries
$$
\begin{aligned}
(STS^{-1})_{1,1}=&\lambda=(STS^{-1})_{2,2}^{-1} \\
(STS^{-1})_{2,1}=&(STS^{-1})_{2,3}=(STS^{-1})_{3,1}=0\\
(STS^{-1})_{1,2}=&\lambda\mu \beta+\mu^2\alpha-\mu^2A_2^+B_4^{-1}B_3^-+\lambda^{-1}\mu \beta-
                 \mu B_2^+A_3^--\mu B_2^+A_4B_4^{-1}B_3^- \\
(STS^{-1})_{1,3}=&-\lambda B_2^+B_4^{-1}+\mu A_2^+B_4^{-1}+B_2^+A_4B_4^{-1}\\
(STS^{-1})_{3,2}=&\lambda^{-1}\mu B_3+\mu B_4^{-1}A_3^--\mu B_4A_4B_4^{-1}B_3^- \\
(STS^{-1})_{3,3}=&B_4A_4B_4^{-1}
\end{aligned}
$$  
In particular, if $|\lambda|\neq1$, by choosing $\mu=1$, $B_4=\id$, 
there exists $B_3^-$ such that
$$A_3^-+(\lambda^{-1}-A_4)(B_3^-)=0\,.$$
Then, using the relations (i) and (ii), one can see that $T$ is conjugate to 
$$
\bpm
\lambda&0&0\\
0&\lambda^{-1}&0\\
0&0&A_4
\epm\,.
$$
and is hence hyperbolic; conversely, one can show that if
$|\lambda|\neq1$, $T$ is hyperbolic.

%

\end{document}